%% file: w-types-std.tex
\documentclass{lmcs}
\pdfoutput=1

\usepackage{lastpage}
\lmcsdoi{17}{3}{28}
\lmcsheading{}{\pageref{LastPage}}{}{}%
{Sep.~16,~2019}{Sep.~24,~2021}{}

\usepackage[utf8]{inputenc}

\usepackage{microtype}
\usepackage{amsmath,amssymb,amsthm}
\usepackage{prooftree}
\usepackage[all]{xy}
\usepackage{mathtools}
\usepackage{xspace}

\usepackage{hyperref}
\usepackage[capitalise,noabbrev]{cleveref}

\input{w-types-std-macros.tex}

\begin{document}

\title{W-types in setoids}

\author{Jacopo Emmenegger}
\address{Matematiska institutionen, Stockholms unversitet, Sweden.}
\curraddr{DIMA, Università degli Studi di Genova, 16146 Genova, Italy.}
\email{emmenegger@dima.unige.it}
\thanks{This work was partially funded by EPSRC grant EP/T000252/1.}

\begin{abstract}
We present a construction of W-types in the setoid model
of extensional \mltt
using dependent W-types in the underlying intensional theory.
More precisely, we prove that the internal category of setoids
has initial algebras for polynomial endofunctors.
In particular, we characterise the setoid of algebra morphisms
from the initial algebra to a given algebra
as a setoid on a dependent W-type.
We conclude by discussing the case of free setoids.
We work in a fully intensional theory and,
in fact, we assume identity types only when discussing free setoids.
By using dependent W-types
we can also avoid elimination into a type universe.
The results have been verified in Coq and a formalisation is available on
the author's GitHub page.
\end{abstract}

\maketitle

\section{Introduction}

The present paper is a contribution to the study of
models of extensional properties in intensional type theories
and is in particular concerned with W-types.
The W-type constructor in \mltt~\cite{Martin-Loef1984,Martin-Loef1982}
produces an inductive type
whose terms can be understood as well-founded trees.
We provide a construction, verified in Coq,
of W-types in the setoid model
of extensional \mltt using dependent W-types
in the underlying intensional theory.
Although we work internally in intensional \mltt,
we present our results using the category-theoretic language.
More precisely, we consider a logical-framework presentation
of a dependent type theory with unit type, $\sum$-types, $\prod$-types
and dependent W-types and we show in \cref{thm:init} that
the internal category of setoids
has initial algebras for polynomial endofunctors.
These were identified as a category-theoretic counterpart of W-types
by Moerdijk and Palmgren~\cite{MoerdijkPalmgren2000}.
Initiality of an algebra and the induction principle of the underlying type
are related via the notion of contractibility in~\cite{AGS2017}.

Dependent W-type were introduced by Petersson and Synek~\cite{PeterssonSynek1989},
see also~\cite[Chapter~16]{NPS1990},
to provide a constructor for general inductive data types,
and are also known as indexed W-types or tree types.
The dependent W-type constructor produces a family of mutually inductive types,
as opposed to the single inductive type produced by the ordinary W-type constructor.
Indeed, ordinary W-types are type-theoretically equivalent to dependent W-types
indexed over the unit type.
Gambino and Hyland~\cite{GambinoHyland2004}, see also~\cite{GambinoKock2013},
formulated a dependent version of polynomial endofunctors
and identified in their initial algebras
a category-theoretic counterpart of dependent W-types.

By setoid we mean a type equipped
with a type-valued total equivalence relation.
That these form a model of extensional \mltt that
interprets most type constructors
has been known for some time.
Most type constructors, as only recently Palmgren~\cite{Palmgren2019}
has provided us with a solution to the
long standing problem of interpreting a type universe,
and in fact a whole hierarchy,
in the setoid model.
As to W-types, it was again Palmgren who first showed
how to construct W-types for setoids
in intensional \mltt, see~\cite{Bressan2015}.
His approach uses extensively the elimination principle of W-types
into a type universe.
The novelty of our approach consists
in the use of dependent W-types instead.

We use dependent W-types in two crucial steps.
First, to define the (partial) equivalence relation on a W-type
that gives rise to the initial algebra setoid $W$.
Second, to characterise the setoid of morphisms of algebras
from the (candidate) initial algebra to a given algebra.
After these steps, it only takes an easy induction
to conclude initiality of $W$.
To briefly describe the steps where dependent W-types are used,
let us recall that we may regard terms of a W-type as well-founded trees:
a canonical tree is recursively specified providing its root and its immediate subtrees,
\ie other well-founded trees that are to be connected to the root.

The relation that we define in \ref{def:Wper}
stipulates that two trees are
matching if their roots are equal and
their immediate subtrees on equal branches are matching,
where by equal we mean the setoid equality.
By using dependent W-types,
this inductively defines a partial relation on the underlying W-type,
and we say that a tree is \emph{extensional}
if it is matching with itself.
Palmgren's construction of W-types for setoids,
as well as the other constructions that we discuss below,
use the same relation, but constructed in different ways.
The initial algebra will be the setoid of extensional trees
with algebra map given by the constructor of the W-type.

Given an algebra $a$ on a setoid $A$,
we then need to construct an algebra morphism $f$ from $W$ to $A$,
that is, a function commuting with the algebra maps,
and to show that it is the unique such.
Commutativity tells us, roughly, that
the action of $f$ on an extensional tree $w$ is determined,
via the algebra $a$, by the action of $f$ on the immediate subtrees of $w$.
By using dependent W-types, we make precise in \cref{def:tel} what it means
for a function on extensional trees
to be determined by its action on immediate subtrees,
and call these functions \emph{telescopic}.
Actually, for the construction itself we find it more convenient
to speak of telescopic functions on immediate subtrees rather than on trees,
but here we can safely ignore the difference.
In \cref{thm:commchar} we prove that a function on extensional trees
is an algebra morphism \iff it is telescopic.

The characterisation in \cref{thm:commchar} allows us to reduce the problem of finding
a unique algebra morphism, to the problem of finding a unique telescopic function.
To this aim, and thanks to the inductive nature of telescopic functions,
we directly use the elimination principle of dependent W-types.
In this sense, we believe that \cref{thm:commchar} makes explicit,
in the case of setoids, the connection between the commutativity condition
for an algebra morphism out of the initial algebra,
and the inductive definition of the morphism itself.

In particular, a common aspect of arguments that construct W-types
is the use of the set, or setoid, of all subtrees of a tree,
usually obtained as the transitive closure of the immediate subtree relation.
This set may have a complicated construction in categories
and in intensional type theories,
but it is important for inductive arguments.
We can avoid dealing with transitive closures when proving \cref{thm:commchar}
since the definition of telescopic functions
only involves the setoid of immediate subtrees of a tree,
which has a natural and straightforward
definition as an image factorisation,
see \cref{def:imsub,rem:imgfact}.

We also use dependent W-types to compare
the setoid of extensional trees \wrt discrete setoids
to the discrete setoid on the W-type of the underlying types.
In particular, we see in \cref{cor:hset}
that every discrete tree with nodes from a 0-type,
\ie a type with decidable equality,
is extensional.
It seems that, without further assumptions,
not every discrete tree is extensional.
However, if function extensionality holds,
then it is possible to identify
the setoid of extensional trees on discrete setoids
with a subsetoid of discrete trees, see~\cref{thm:freestd}.

The setoid construction that we consider is an instance
of a quotient completion,
see~\cite{MaiettiRosolini2013,MaiettiRosolini2016}.
The author is aware of two other constructions of W-types
for quotient completions.
The first one, that we already mentioned,
was formulated by Palmgren for the setoid model in intensional \mltt
and then adapted by Bressan~\cite{Bressan2015} to
minimalist type theory~\cite{Maietti2009}.
The argument requires a `large' elimination principle for W-types,
in the sense that it must be possible to eliminate into a type universe
or a universe of propositions.
The second construction is due to van den Berg~\cite{vdBerg2005}
and it applies to exact completions of categories with finite limits.
In intensional \mltt, the assumption of finite limits is met
by the category of types
assuming function extensionality and Uniqueness of Identity Proofs,
by its full subcategory on the 0-types
assuming only function extensionality~\cite{HoTTbook,RijkeSpitters2015},
and by the e-category of types assuming only UIP~\cite{EmmeneggerPalmgren2017}.
However, there is little hope that an internal category of types in a
fully intensional type theory would have ordinary finite limits.
Setoids are also closely related to groupoids.
Vidmar has given a construction of initial algebras
for polynomial endofunctors on groupoids
from initial algebras for polynomial endofunctors on sets~\cite{Vidmar2018}.

In the next section we discuss the preliminaries needed for the construction.
\Cref{sec:exttree} contains the construction of the algebra of
extensional trees and the definition of the setoid family of
immediate subtrees.
The proof of its initiality is in \cref{sec:init}.
We conclude the paper with a discussion of the case of
extensional trees over discrete setoids in \Cref{sec:disctree}.
Each section begins with a brief overview of the content.

The Coq formalisation of the present paper is available on the author's GitHub page~\cite{Emmenegger2018a} and it includes
\cref{def:std,def:extfun} and
all numbered definitions and results in
\cref{ssec:poly,sec:exttree,sec:init,sec:disctree},
except for \cref{thm:lccstd},
a Coq proof of which can be found at~\cite{Palmgren2012c}.

\section{Preliminaries}
\label{sec:prelim}

We present the type theory we will be working with in \cref{ssec:ttbasics};
briefly review W-types in \cref{ssec:wty}
and their dependent verion in \cref{ssec:dwty};
and recall some basic facts and definitions
about setoid and setoid families in \cref{ssec:stdfam},
and about polynomial functors in \cref{ssec:poly}.

\begin{figure}[b]
\caption{Rules for W-types.}
\label{wrules}
\begin{gather*}
\begin{prooftree}
A \typing \ttunivm	\qquad	B \typing A \to \ttunivm
\justifies \rule{0pt}{13pt}
\wty[A,B] \typing \ttunivm
\using{\footnotesize \textsf{W-FORM}}
\end{prooftree}
\qquad \qquad
\begin{prooftree}
a \typing A	\qquad	f \typing B \fapp a \to \wty[A,B]
\justifies \rule{0pt}{13pt}
\wsup_{A,B} \fapp a \fapp f \typing \wty[A,B]
\using{\footnotesize \textsf{W-INTRO}}
\end{prooftree}
\\[3ex]
\begin{prooftree}
C \typing \wty[A,B] \to \ttunivm	\qquad
c \typing \prod_{(a \typing A)} \prod_{(f \typing B \fapp a \to \wty[A,B])}%
    \Bigg( \prod_{b \typing B \fapp a}%
	    C\fapp (f \fapp b) \Bigg) \to C \fapp (\wsup_{A,B} \fapp a \fapp f)
\justifies \rule{0pt}{13pt}
\wrec[A,B,C,c] \typing \prod_{w \typing \wty[A,B]} C \fapp w
\using{\footnotesize \textsf{W-ELIM}}
\end{prooftree}
\\[3ex]
\begin{prooftree}
a \typing A	\qquad	f \typing B \fapp a \to \wty[A,B]
\justifies \rule{0pt}{13pt}
\wrec[A,B,C,c] \fapp (\wsup_{A,B} \fapp a \fapp f)%
    \jueq c \fapp a \fapp f \fapp (\lambda b.\wrec[A,B,C,c] \fapp (f \fapp b))
\using{\footnotesize \textsf{W-CONV}}
\end{prooftree}
\end{gather*}
\end{figure}

\begin{figure}
\caption{Rules for dependent W-types.}
\label{dwrules}
\begin{gather*}
\begin{prooftree}
\[ I \typing \ttunivm	\qquad	N \typing I \to \ttunivm
\justifies \thickness=0pt
Br \typing \prod_{i \typing I} N \fapp i \to \ttunivm	\qquad%
ar \typing \prod_{(i \typing I)}
\prod_{(n \typing N \fapp i)} Br(i,n) \to I \]
\justifies \rule{0pt}{13pt}
\dwty[Br,ar] \typing I \to \ttunivm
\using{\footnotesize \textsf{DW-FORM}}
\end{prooftree}
\\[3ex]
\begin{prooftree}
i \typing I	\qquad	n \typing N \fapp i	\qquad
f \typing \prod_{b \typing Br(i,n)} \dwty[Br,ar] \fapp ar(i,n,b)
\justifies \rule{0pt}{13pt}
\wsupd_{Br,ar}(i,n,f) \typing \dwty[Br,ar] \fapp i
\using{\footnotesize \textsf{DW-INTRO}}
\end{prooftree}
\\[3ex]
\resizebox{\columnwidth}{!}{%
\begin{prooftree}
\[ C \typing \prod_{i \typing I} \dwty[Br,ar] \fapp i \to \ttunivm
\hspace{.777\textwidth}\justifies \thickness=0pt
c \typing \prod_{(i \typing I)} \prod_{(n \typing N \fapp i)}%
    \prod_{\left(f \typing \prod_{b \typing Br(i,n)}%
	\dwty[Br,ar] \fapp ar(i,n,b)\right)}%
\Bigg( \prod_{b \typing Br(i,n)} C(ar(i,n,b),f \fapp b) \Bigg)%
    \to C(i,\wsupd_{Br,ar}(i,n,f)) \]
\justifies \rule{0pt}{13pt}
\dwrec[Br,ar,C,c] \typing \prod_{i \typing I} \prod_{w \typing \dwty[Br,ar] \fapp i}%
    C(i,w)
\using{\footnotesize \textsf{DW-ELIM}}
\end{prooftree}
}
\\[3ex]
\begin{prooftree}
i \typing I	\qquad	n \typing N \fapp i	\qquad
f \typing \prod_{b \typing Br(i,n)} \dwty[Br,ar] \fapp (ar(i,n,b))
\justifies \rule{0pt}{13pt}
\dwrec[Br,ar,C,c](i,\wsupd_{Br,ar}(i,n,f))%
    \jueq c(i,n,f,\lambda b.\dwrec[C,c](ar(i,n,b) ,f \fapp b))
\using{\footnotesize \textsf{DW-CONV}}
\end{prooftree}
\end{gather*}
\end{figure}

\subsection{The type-theoretic setting}
\label{ssec:ttbasics}

We work in a logical framework formulation of
\mltt, see~\cite[Part III]{NPS1990},
similar to the fragment of Coq that we used for the formalisation.
The logical framework consists of record types,
$\prod$-types and a universe \`a la Russell $\ttunivm$,
which is the type $\mathsf{Set}$ in the Coq formalisation.
Logic is interpreted according to the propositions-as-types interpretation.
When we wish to emphasise that a type is in $\ttunivm$
we say that it is a \emph{small type}.
Judgemental equality is denoted $\jueq$ and
application of function terms will usually be denoted
by concatenation as in $f\fapp b$.
However, in the presence of multiple arguments we may adopt the notation with parenthesis
as in $ar(i,n,b)$.

The theory itself is then specified by declaring constants and equations.
Specifically, we require the universe \ttunivt
to be closed under $\prod$-types,
and to contain the unit type \ttunitt, $\sum$-types,
and dependent W-types.
Ordinary W-types are recovered as dependent W-types over the unit type  \ttunitt,
but in the Coq formalisation these are assumed to be primitive
for convenience.
In \cref{wrules,dwrules} we spell out
the constants and equations for (dependent) W-types
in the form of rules to increase readability.
In these rules, the premises of \textsf{W-FORM} and \textsf{DW-FORM}
are left understood in the other rules.
We might as well drop (some) subscripts if they are clear from the context.
The system has $\eta$-conversion only for $\prod$-types.

Note that we do not require identity types.
These will be assumed only in \cref{sec:disctree}
for discussing the case of discrete setoids.
It is well-known that record types can be replaced with $\sum$-types
at the expense of readability.
In this case we would not need to assume that these coincide
with the $\sum$-types in the theory.

\subsection{W-types}
\label{ssec:wty}

W-types can be used to construct several inductive types,
including natural numbers and lists~\cite{Dybjer1997},
and to give a constructive justifications
of certain theories of inductive definitions~\cite{Palmgren1992}.
More generally, W-types provide a predicative counterpart to the notion of well-ordering. 
Furthermore, they are instrumental in Aczel's model of
Constructive Zermelo-Fraenkel set theory~\cite{Aczel1978},
in the form of the type of iterative sets,
where they also allow to interpret the Regular Extension Axiom
which adds general inductive definitions to CZF~\cite{Aczel1986}.
The idea of a type of iterative sets is also central
in Palmgren's interpretation of type universes in the setoid model~\cite{Palmgren2019}.

By elimination of W-types, there are function terms
\[
\wn \typing \wty[B] \to A
\qquad\text{and}\qquad
\wb \typing \prod_{w \typing \wty[B]} B \fapp (\wn \fapp w) \to \wty[B]
\]
such that $\wn \fapp (\wsup \fapp a \fapp f) \jueq a $
and $\wb \fapp (\wsup \fapp a \fapp f) \jueq f$.
Henceforth, we write $\wb_w$ for $\wb\fapp w$.

We refer to terms of a W-type $\wty[B]$ as \emph{trees},
to terms $a \typing A$  as \emph{nodes},
to terms $b \typing B\fapp a$  as \emph{branches} and,
for every branch $b \typing B\fapp a$, to the term $\wb_w\fapp b$ as the
\emph{immediate subtree of $w$ on branch $b$}.

In the set-theoretic interpretation of type theory,
W-types are sets of well-founded trees with labelled nodes.
The set $A$ is the set of labels for the nodes
and, for each $a \in A$,
the set $B_a$ consists of the branches out of the node with name $a$.
Trees are formed by providing a node $a \in A$
and attaching other trees to the branches in $B_a$.
This procedure is formally specified by functions $f \colon B_a \to \wty[B]$
which stipulate that the tree $f(b)$ is attached to the branch $b \in B_a$.

\begin{figure}[t]
\caption{Canonical term of a W-type.}
\label{tree}
\[\xycenterm{
&	\bullet	\ar@{-}[dr]_<{a_0}
&&	\bullet	\ar@{-}[dl]^<{a_0}
\\
\bullet	\ar@{-}[dr]_<{a_0}
&&	\bullet	\ar@{-}[dl]^<{a_2}^>{a_1}
&\\
&	\bullet
&&\\
}\hspace{6em}
\wsup\fapp a_1\fapp \left(\begin{split}
\strut
0 &\mapsto \wsup\fapp a_0\fapp f_{\emptyset}
\\
1 &\mapsto \wsup\fapp a_2\fapp \left(\begin{split}
0 &\mapsto \wsup\fapp a_0\fapp f_{\emptyset}
\\
1 &\mapsto \wsup\fapp a_0\fapp f_{\emptyset}
\end{split}\right)
\end{split}\right)
\]
\end{figure}

\begin{figure}[b]
\caption{Canonical term of a dependent W-type.}
\label{deptree}
\[\xycenterm[R=1em]{
&	\bullet	\ar@{-}[d]_<{a_0}
&&	\bullet	\ar@{-}[d]^<{a_0}
\\
&	*\txt{\strut\scriptsize $\,i_0\,$}	\ar@{-}[dr]
&&	*\txt{\strut\scriptsize $\,i_0\,$}	\ar@{-}[dl]_>{a_2}
\\
\bullet	\ar@{-}[d]_<{a_0}
&&	\bullet	\ar@{-}[d]
&\\
\txt{\scriptsize $\,i_0\,$}	\ar@{-}[dr]
&&	\txt{\scriptsize $\,i_1\,$}	\ar@{-}[dl]_>{a_1}
&\\
&	\bullet	\ar@{-}[d]
&&\\
&	\txt{\scriptsize $\,i_1\,$}	&&
}\hspace{8ex}
\wsupd\fapp i_1\fapp a_1\fapp \left(\begin{split}
\strut
0 &\mapsto \wsupd\fapp i_0\fapp a_0\fapp f_{\emptyset}
\\
1 &\mapsto \wsupd\fapp i_1\fapp a_2\fapp \left(\begin{split}
0 &\mapsto \wsupd\fapp i_0\fapp a_0\fapp f_{\emptyset}
\\
1 &\mapsto \wsupd\fapp i_0\fapp a_0\fapp f_{\emptyset}
\end{split}\right)
\end{split}\right)
\]
\end{figure}

\Cref{tree} shows a tree in $\wty[B]$ when
$A$ contains three elements $a_0$, $a_1$ and $a_2$
such that $B_{a_0}$ is empty and $B_{a_1} = \{0,1\} = B_{a_2}$.
The corresponding canonical element of $\wty[B]$ is displayed on the right-hand side
together with the five functions mapping a branch to the tree attached to it,
among which occur the empty function $f_{\emptyset}$ into $\wty[B]$.

As it is well-known, these sets of trees
can be characterised more abstractly
as free term algebras for infinitary single-sorted signatures.
In this view,
the set $A$ contains the function symbols of the signature,
while (the cardinality of) $B \fapp a$ is the arity of symbol $a$.
Terms are built out of function symbols according to
composition instructions specified by functions
$f \colon B_a \to \wty[B]$
which stipulate that, for $b \in B_a$,
the term $f(b)$ occurs as $b$-th argument.
The tree in \cref{tree} is the term also written $a_1(a_0,a_2(a_0,a_0))$,
in a signature with function symbols $a_0$ with arity $0$
and $a_1,a_2$ with arity $2$.

\subsection{Dependent W-types}
\label{ssec:dwty}

In the same vein,
dependent W-types can be understood as
free term algebras for infinitary \emph{multi-sorted} signatures:
sorts are elements of $I$,
elements of $N_i$ are function symbols with codomain sort $i $
and the arity of a function symbol $a \in N_i$ is given by the function
$ar_{i,a} \typing B_{i,a} \to I$ that maps each
$b$-th argument of $a$ to its sort.

\Cref{deptree} depicts an element of a dependent W-type.
Compared to the graphical representation of trees in \cref{tree},
the sort $i$ of a function symbol $a\in N_i$ appears as an additional leg
pointing downwards out of the node corresponding to $a$,
and each branch $b \in  Br_{i,a}$ is labelled with the sort $ar(i,a,b) \in I$.
The dependent tree in \cref{deptree} is a term in
the signature containing two sorts $i_0$ and $i_1$,
and three function symbols,
$a_0$ of sort $i_0$ and $a_1,a_2$ of sort $i_1$,
such that $a_0$ has empty arity,
$a_1$ has arity $i_0$ and $i_1$,
and $a_2$ has arity $i_0$ and $i_0$.
This is formally specified by requiring that
$i_0,i_1 \in I$,
$a_0 \in N_{i_0}$, $a_1,a_2 \in N_{i_1}$,
$Br_{i_0,a_0} = \emptyset$,
$Br_{i_1,a_1} = \{0,1\} = Br_{i_1,a_2}$,
$ar(i_0,a_0)$ is the empty function into $I$,
$ar(i_1,a_1)$ maps $0$ to $i_0$ and $1$ to $i_1$ and
$ar(i_1,a_2)$ maps $0$ and $1$ to $i_0$.
A common presentation of the term in \cref{deptree}
is $a_1(a_0,a_2(a_0,a_0)) \typing i_1$.

Note that there must be at least one node with no branches,
\ie at least one constant,
in order for these trees to be well-founded.
Indeed, in terms of signatures, these sets of trees consist of closed terms.
This is reflected in type theory:
whenever $B$ is a family of non-empty types,
$ \wty[B]$ is type-theoretically equivalent to the empty type \ttemptyt, when the latter is available.
In order to make this precise,
the meaning of `non-empty' can be understood either as
$ \prod_{a \typing A} (B \fapp a \to \ttemptym) \to \ttemptym $
or as $\prod_{a \typing A} || B \fapp a ||$
if propositional truncation is available,
in the sense that assuming either of them proves
that $\wty[B]$ is equivalent to $\ttemptym$.

Similarly to the non-dependent case,
using the elimination rule of dependent W-types,
we obtain for every $i \typing I$  two function terms
\begin{gather}
\label{eq:dwn}
\dwn_i \typing \dwty \fapp i \to N \fapp i
\\[1.5ex]
\label{eq:dwb}
\dwb_i \typing \prod_{(w \typing \dwty \fapp i)}
\prod_{(b \typing Br(i,\dwn_i \fapp w))}%
    \dwty \fapp ar(i,\dwn_i \fapp w,b),
\end{gather}
such that $\dwn_i \fapp \wsupd(i,n,f) \jueq n$
and $\dwb_i \fapp \wsupd(i,n,f) \jueq f$.
We refer to $I$ as the type of \emph{indices} or sorts,
to $N$ and $Br$ as the \emph{node family} and the \emph{branching family}, respectively,
and to $ar \typing \prod_i\prod_n Br(i,n) \to I$
as the \emph{arity function}.

The following examples are just for the sake of illustrating the use of
dependent W-types and we do not use the definitions therein elsewhere in the paper.

\begin{xmps}\hfill
\begin{enumerate}
\item
Let $X \typing \ttunivm$ and $Y \typing X \to  \ttunivm$.
We can define a type of (root-to-leaf) paths in a tree $w \typing \wty[Y]$
as a dependent W-type over $\wty[Y]$ as follows.
Say that a path in $w \jueq \wsup\fapp x\fapp f$
is a branch $y \typing Y\fapp x$
together with a path in the immediate subtree
$f\fapp y$.

The family $\mathsf{Path}_Y $ of paths in trees in $\wty[Y]$
is the dependent W-type on $\wty[Y]$ with node and branching families
\[
N \coloneqq \lambda w.Y\fapp (\wn\fapp w) \typing \wty[Y] \to \ttunivm,
\hspace{4em}
Br \coloneqq \lambda w,y.\ttunitm \typing \prod_{w \typing \wty[Y]} N\fapp w \to \ttunivm,
\]
and with arity function
$ar \coloneqq \lambda w,y,\_. \wb_w\fapp y \typing \prod_w\prod_y\prod_{\ttunitm} \wty[Y].$
A canonical element of $\mathsf{Path}_Y \fapp w$ consists of
a branch $y$ of $w$ together with a function
$\ttunitm \to \mathsf{Path}_Y (\wb_w\fapp y)$,
\ie a path in the immediate subtree $\wb_w\fapp y$.

\item
If the theory has finite sums (in fact, the type $\mathsf{2}$ suffices)
and the empty type $\mathsf{0}$,
it is possible to define a type of finite paths in a tree
$w \typing \wty[Y]$
similarly to the previous example.
Say that a finite path in $w \jueq \wsup\fapp x\fapp f$
is either the empty path,
or it is a branch $y \typing Y\fapp x$
together with a finite path in the immediate subtree
$f\fapp y \typing \wty[Y]$.

The family $\mathsf{FinPath}_Y \typing \wty[Y] \to \ttunivm$
of finite paths in trees in $\wty[Y]$ is the dependent W-type
on $\wty[Y]$ with node and branching families
\[
N \coloneqq \lambda w. Y\fapp (\wn\fapp w) + \ttunitm,
\hspace{4em}
Br \coloneqq \lambda w,z.\begin{cases}
\, z \jueq \mathsf{inl}\fapp y &\mapsto \ttunitm
\\
\, z \jueq \mathsf{inr}\fapp t &\mapsto \mathsf{0}
\end{cases}
\]
and with arity function $ar \coloneqq \lambda w,z,\_.\begin{cases}
\, z \jueq \mathsf{inl}\fapp y &\mapsto \wb_w\fapp y
\\
\, z \jueq \mathsf{inr}\fapp t &\mapsto R_\mathsf{0}
\end{cases} \typing \prod_w\prod_z\prod \wty[Y]$
where $R_\mathsf{0} \typing \mathsf{0} \to \wty[Y]$
is the eliminator of the empty type.

A canonical element of $\mathsf{FinPath}_Y \fapp w$ consists
either of a term $t \typing \ttunitm$ and a function
$\mathsf{0} \to  \wty[Y]$,
which we regard as `the' empty path,
or of a branch $b$ of $w$ and a function
$\ttunitm \to \mathsf{FinPaths}_Y\fapp (\wb_w\fapp y)$,
\ie a finite path in the immediate subtree $\wb_w\fapp y$.
Either this path consists of a term of $\ttunitm$
and an empty function, in which case the path is over with length one.
Or it consists of a branch $y'$ and
a finite path in $\wb_{\wb\fapp y}\fapp y'$.

The function term
$\mathsf{FinPath}_Y \fapp w \to \mathsf{N}$
that maps each finite path to its length can be easily defined
using the elimination rule of dependent W-types in \cref{dwrules}.
\end{enumerate}
\end{xmps}

\subsection{Setoids and setoid families}
\label{ssec:stdfam}

A setoid is informally defined as a type together with
a type-theoretic equivalence relation on it.
This can be made precise in various different ways,
\eg considering partial relations,
or Prop-valued relations if a type of proposition is available,
or instead 0-types and mere relations.
See~\cite{BCP2003} for a survey of possible definitions.
In our context we define a setoid as follows.

\begin{dfn}\label{def:std}
A \emph{setoid} $A$ is a list $(|A|, \approx_A, r_A, s_A, t_A)$ where
$ |A| \typing \ttunivm$ is a small type,
the type family $\approx_A \typing |A| \to |A| \to \ttunivm$
is an equivalence relation with proofs of reflexivity, symmetry and transitivity
given by, respectively:
\begin{gather*}
r_A \typing \prod_{a \typing |A|} a \approx_A a,
\hspace{6em}
s_A \typing \prod_{a, a' \typing |A|} a \approx_A a' \to a' \approx_A a,
\\[1.5ex]
t_A \typing \prod_{a, a', a'' \typing |A|}%
	a \approx_A a' \to a' \approx_A a'' \to a \approx_A a''.
\end{gather*}
The \emph{type of setoids} \catstdt is defined as a record on the types of
$ |A|, \approx_A, r_A, s_A$ and $t_A$.

For a setoid $A$, we refer to the relation $\approx_A$
as the \emph{equality of $A$}.
\end{dfn}

Since our theory has  a unit type,
every small type $X$ gives rise to a setoid on $X$
with equality $\lambda a,a'.\ttunitm$,
called the the \emph{codiscrete setoid on $X$}.

In \cref{sec:disctree} we will have identity types available.
In this case every small type $X$ can be equipped with the equality
$\lambda x,x'. x =_X x'$ given by the identity type.
This setoid is called the \emph{discrete setoid on $X$}.
More generally, a setoid $A$ is \emph{discrete} or (\emph{free}) if
$a \approx_A a'$ \iff $a =_{|A|} a'$
for every $a,a' \typing A$.

In the category of setoids that we will define below,
the codiscrete setoid over any inhabited type is a terminal object,
and it is isomorphic to the discrete setoid on \ttunitt
if identity types are available.

\begin{dfn}\label{def:extfun}
Let $A$ and $B$ be setoids.
A function term $f \typing |A| \to |B|$ is
\emph{extensional (\wrt $\approx_A$ and $\approx_B $)}
if there is a term of type
\begin{equation}\label{eq:extfun}
\mathsf{ext}\fapp f \coloneqq \prod_{a,a'\typing |A|}
a \approx_A a' \longrightarrow f\fapp a \approx_B f\fapp a'.
\end{equation}
We will refer to extensional function terms as extensional functions,
or simply as functions, and
we will write $|f|$ for the underlying function term
of an extensional function $f$.

The setoid $A \extfnc B$ of \emph{extensional functions} from $A$ to $B$
has the type of extensional functions
\[
\sum_{f \typing |A| \to |B|} \mathsf{ext}\fapp f
\]
as underlying type and the equivalence relation
\[
f \approx_{A \extfnc B} g \coloneqq \prod_{a \typing |A|} |f|\fapp a \approx_B |g|\fapp a
\]
as equality.
\end{dfn}

In fact, in the Coq implementation we find it more convenient to define the type
of extensional functions as a record rather than as a $\sum$-type.
However this makes no essential difference.

\begin{rmk}\label{rem:extper}
Both type families $\mathsf{ext}$ and $\approx_{A \extfnc B}$
defined in \ref{def:extfun} are instances of a more general
binary relation on the type $|A| \to |B|$, namely
\[
\mathsf{extgen}(f,g) \coloneqq \prod_{a,a' \typing |A|}
a \approx_A a' \longrightarrow f\fapp a \approx_B g\fapp a'.
\]
It is easy to see that this is a partial equivalence relation on $|A| \to |B|$
whose domain, \ie the type of function terms $f$ such that $\mathsf{extgen}(f,f)$,
is the type of extensional function $|A \extfnc B|$.
Furthermore, two extensional functions $f$ and $g$ are equal
\iff the type $\mathsf{extgen}(|f|,|g|)$ is inhabited.
\end{rmk}

In the rest of the paper we write $a \typing A$ for a setoid $A$,
to mean $a \typing |A|$,
and we often not distinguish between
an extensional function $f \typing A \extfnc B$
and its underlying function term $|f| \typing |A| \to |B|$.
We also write $f\extpf \alpha \typing f \fapp a \approx_B f \fapp a' $
for the proof of extensionality of $f$ applied to $\alpha \typing a \approx_A a'$.
Occasionally, we also find it convenient
to drop the subscript from the equality of a setoid.
We do not expect these abuses of notation to lead to confusion.

Setoids provide us with a formulation of the notion of category 
without equality on objects that was introduced by Aczel,
who formalised it in Lego~\cite{Aczel1993}.
See also~\cite{HuetSaibi2000} for a formalisation in Coq
that uses setoids with Prop-valued equalities.

\begin{dfn}\label{def:ecat}
A \emph{(locally small) e-category} \catat consists of
a type of objects $\Obj_\catam$
and a family of setoids
$\Hom_\catam \typing \Obj_\catam \to \Obj_\catam \to \catstdm$,
of arrows,
together with function terms for identity and composition,
where the latter is extensional, in the sense that it has type
\[
\prod_{a,b,c \typing \Obj_\catam} \Hom_\catam(b,c) \extfnc \Hom_\catam(a,b) \extfnc \Hom_\catam(a,c),
\]
and with identity and associativity axioms formulated
using equalities of the family of setoids $\Hom_\catam$.

An \emph{e-functor} between two e-categories \catat and \catbt consists of
a function term $F \typing \Obj_\catam \to \Obj_\catbm$ 
together with a term of type
\[
\prod_{a,a' \typing \Obj_\catam} \Hom_\catam(a,a') \extfnc \Hom_\catbm(F \fapp a, F \fapp a')
\]
and proof terms witnessing that it is functorial
with respect to the equalities of the family of setoids $\Hom_\catbm$.

An \emph{e-natural transformation} between two e-functors $F,G \colon \catam \to \catbm$ 
consists of a term
\[
n \typing \prod_{a \typing \Obj_\catam} \Hom_\catbm(F \fapp a, G \fapp a)
\]
whose action on $a \typing \Obj_\catam$ we denote as $n_a$, together with a term of type
\[
\prod_{a,a' \typing \Obj_\catam}\prod_{f \typing a \extfnc a'}
G \fapp  f \circ n_a \approx n_{a'} \circ F \fapp f,
\]
where $\approx$ is the equality of $\Hom_\catbm(F \fapp a,G \fapp a')$.
\end{dfn}

The names `e-category', `e-functor', etc.\  were introduced to distinguish
these concepts from the standard ones.
However, since these are the only formulations that appear in the present paper,
we will just say `category' to mean `e-category' and similarly for the others.

For a setoid $A$, the \emph{discrete category} $A^\#$ on $A$ is defined as follows.
Its type of objects is $|A|$ and, for $a, a' \typing A$,
its setoid of arrows from $a$ to $a'$ is the codiscrete setoid on
the type $a \approx_A a'$.
Since $\approx_A$ is symmetric, the category $A^\#$ is a groupoid,
\ie a category where every arrow is invertible.

Functors between two categories \catat and \catbt
and natural transformations between them
form a category $\efun{\catam,\catbm}$.
We denote the action of a functor $F$ on an arrow $\alpha \typing \Hom_\catam(a,a')$
as $F_\alpha \typing \Hom_\catbm(F \fapp a, F \fapp a')$.

Setoids and extensional functions form a category~\cite{PalmgrenWilander2014}.

\begin{dfn}\label{def:catstd}
The \emph{category of setoids} is defined as follows.
The type of objects is the type of setoids \catstdt
and the family of arrows is the family $\lambda A,B. A \extfnc B$
of setoids of extensional functions.
Identity and composition are defined in the obvious way,
the latter will be denoted as $g \fcmp f$.
We abuse notation and denoted this category also as \catstdt.

For a setoid $A$, we define
\[
\sfam A \coloneqq \efun{A^\#,\catstdm}
\]
and call a functor $B \typing \sfam A$  a \emph{setoid family over $A$}.
We refer to the action of $B$
on an arrow $\alpha \typing a \approx_A a'$ of $A^\#$
as \emph{transport along $\alpha$}.
For $b \typing B\fapp a$ and $b' \typing B\fapp a'$ we define
\begin{equation}\label{def:catstd:depeq}
b \approx_{\alpha} b' \ \coloneqq\
B_{\alpha} \fapp b \approx_{B \fapp a'} b'.
\end{equation}

For every extensional function $f \typing A' \extfnc A$
and every setoid family $B \typing \sfam A$,
there is a setoid family $f^*B \typing \sfam A'$ defined by
$(f^*B)\fapp a \coloneqq B\fapp (f\fapp a)$ and
$(f^*B)_{\alpha} \coloneqq B_{(f\extpf \alpha)}$.
This action extends to a functor from $\sfam A$ to $\sfam A'$.
\end{dfn}

\begin{rmk}\label{rem:famirr}
Besides the standard functoriality conditions on identities and composites,
setoid families enjoy an additional coherence.
Indeed, for every $\alpha,\alpha' \typing a \approx_A a'$,
it is $\alpha \approx \alpha'$  as arrows in $A^\#$.
It follows by extensionality of $B$ on arrows that
\[
B_{\alpha} \approx_{(B \fapp a \extfnc B \fapp a')} B_{\alpha'}
\]
for every $\alpha$  and $\alpha'$  as above.
Indeed, setoid families in this sense are
the proof-irrelevant setoid families studied in~\cite{Palmgren2012a}.
In particular, for every $\alpha, \alpha' \typing a \approx_A a'$,
$ b \typing B\fapp a$  and $b' \typing B\fapp a'$, it is
\begin{equation}\label{eq:famirr}
b \approx_{\alpha} b' \longleftrightarrow b \approx_{\alpha'} b'
\end{equation}
and, for every $\alpha \typing a \approx_A a$ and $b,b' \typing B\fapp a$, it is
\begin{equation}\label{eq:famirrref}
b \approx_{B\fapp a} b' \longleftrightarrow b \approx_{\alpha} b'.
\end{equation}

In fact, setoid families are formalised in~\cite{Emmenegger2018a}
directly as a record on the types of function terms and their functoriality
and proof-irrelevance properties.
\end{rmk}

It is an important feature of the category of sets in set theory
that the slice over a set $A$ is equivalent to families of sets indexed over $A$.
The same holds for the category of setoids in the fragment of intensional \mltt that we are considering, see~\cref{ssec:ttbasics}.

\begin{thr}\label{thm:stdfam}
For a setoid $A$, the slice category $\catstdm/A$ 
is equivalent to the category of setoid families $\sfam A$.
\end{thr}

\begin{proof}
We describe the two functors,
leaving the straightforward verification of the equivalence to the reader.

Every setoid family $B \typing \sfam A$ gives rise to a setoid
$(\sum_{a \typing A} B \fapp a, \approx)$
where the equality is
\[
(a,b) \approx (a',b') \coloneqq \sum_{\alpha \typing a \approx_A a'}%
	b \approx_{\alpha} b'.
\]
The projection into $A$ is extensional,
with extensionality proof also given by the first projection.
This defines an extensional function into $A$.

A natural transformation $n$ from $B$ to $C$ is mapped
to the extensional function with underlying term
$ \lambda (a,b).(a,n_a \fapp b)$.
Its extensionality follows,
for every $(\alpha,\beta) \typing (a,b) \approx (a',b')$, from
\[
C_\alpha (n_a \fapp b) \approx_{C \fapp a'} n_{a'} \fapp (B_\alpha b)
\approx_{C \fapp a'} n_{a'} \fapp b'
\]
which holds by naturality of $n$  and extensionality of $n_{a'}$.

Conversely, let $f \typing B \extfnc A$ be a function.
Consider, for every $a \typing A$, the subsetoid of $B$
on those $b$ such that $f\fapp b \approx_A a$, \ie the type
\[
\sum_{b \typing B} f \fapp b \approx_A a
\]
with equality
$(b,\_) \approx (b',\_) \coloneqq b \approx_B b'$.
This assignment extends to a setoid family
whose transport along $\alpha \typing a \approx_A a'$ maps
a pair $(b,\beta)$
to the pair $(b,\beta')$, where $\beta' \typing f \fapp b \approx_A a'$ 
is the concatenation of $\beta$  and $\alpha$.

An arrow in $\catstdm/A$ from $f$ to some $g \typing C \extfnc A$ 
is a function $k \typing B \extfnc C$ such that $g \circ k \approx f$.
This gives rise to a natural transformation as follows.
For $a \typing A$, the underlying function term maps
$(b,\beta)$ to $(k\fapp b,\beta')$,
where $\beta' \typing g \fapp (k \fapp b) \approx_A a$ 
is the concatenation of $\beta \typing f \fapp b \approx_A a$ 
with the proof of $g \circ k \approx f$  applied to $b$.
This function term is clearly extensional as $k$ is.
As transport along $\alpha \typing a \approx_A a'$ 
does not modify the first component,
naturality in $a \typing A$ is trivial.
\end{proof}

\subsection{Polynomial functors and W-types}
\label{ssec:poly}

Recall that a category is \emph{locally cartesian closed} if
it has a terminal object and
all its slices are cartesian closed, that is,
have all binary products and exponentials.
\Cref{thm:stdfam} may be used to prove the following result,
another proof can be found in~\cite{EmmeneggerPalmgren2017}.

\begin{prp}\label{thm:lccstd}
The category of setoids \catstdt is locally cartesian closed.
\end{prp}

In a locally cartesian closed category \catct, for every arrow $f \colon B \to A$,
the pullback functor $f^*$ has a left adjoint $f_!$,
defined by post-composition with $f$,
and a right adjoint $f_*$,
as below.
\begin{equation}\label{eq:lccadj}
\xycenterm[C=6em]{
\catcm/B	\ar[r];[]|-{\,f^*\,} \ar@<-1.5ex>[r]_-{f_!} \ar@<1.5ex>[r]^-{f_*}
&	\catcm/A
}\end{equation}

When $f$ is the unique arrow $B \to 1$ to the terminal object,
the above functors become
\[\xycenterm[C=6em]{
\catcm/B	\ar[r];[]|-{\ (-) \times B\ } \ar@<-1.5ex>[r]_-{B_!} \ar@<1.5ex>[r]^-{(-)^B}
&	\catcm.
}\]
It follows that every arrow $f \colon B \to A$ in \catct
gives rise to a functor $P_f \colon \catcm \to \catcm$,
the \emph{polynomial endofunctor associated to $f$},
defined as the composite
\begin{equation}\label{eq:poly}
\xycenterm[C=4em]{
\catcm	\ar[r]^-{(-)\times B}	&	\catcm/B	\ar[r]^-{f_*}
&	\catcm/A	\ar[r]^-{A_!}	&	\catcm.
}\end{equation}
A \emph{polynomial endofunctor} is then defined to be an endofunctor on \catct
which is naturally isomorphic to $P_f$ for some $f$ in \catct.

An \emph{algebra for a polynomial endofunctor $P$}
is given by an object $X$ and
an arrow $s \colon P X \to X$, called \emph{algebra map}.
Such an algebra is \emph{initial} if for any algebra $t \colon P Y \to Y$,
there is $h \colon X \to Y$ such that $t \fcmp (P \fapp h) = h \fcmp s$
and $h$ is the unique such.
It is a well-known result by Lambek that the algebra map of an initial algebra is invertible.

In extensional type theory with one universe,
the category of small types and function terms is locally cartesian closed
if the universe has \ttunitt, identity types, $\sum$ and $\prod$ types.
Hence we may consider polynomial endofuctors
for each function term $f \typing B \to A$.
An initial algebra is then given by the W-type of the type family
$ f^{-1} \coloneqq \lambda a. \sum_{b \typing B} f(b) =_A a$,
with algebra map defined by
\[
(a, k) \typing \sum_{a \typing A} \Big( f^{-1}(a) \to \wty[f^{-1}] \Big)%
    \longmapsto \wsup(a,k) \typing \wty[f^{-1}].
\]
It is not difficult to see that initiality of $(\wty[f^{-1}], \wsup)$
is logically equivalent to
the recursion principle of $\wty[f^{-1}]$~\cite{MoerdijkPalmgren2000}.

In \catstdt, we take advantage of \cref{thm:stdfam}
and formulate the notion of polynomial functor for setoid families.

\begin{dfn}\label{def:poly}
Let $B \typing \sfam A$ be a setoid family over a setoid $A$.
The \emph{polynomial endofunctor associated to $B $}
is the functor $\catstdm \to \catstdm$ 
defined on $X \typing \catstdm$ by
\[
\polysub{B} X \coloneqq \Big( \sum_{a \typing A} (B \fapp a \extfnc X),\, \approx_{\polysub{B} X} \Big),%
\quad \text{where}\quad
(a, k) \approx_{\polysub{B} X} (a', k') \coloneqq  \sum_{\alpha \typing a \approx a'}%
  k \approx k' \fcmp B_\alpha
\]
and on $f\typing X \extfnc Y$ by
\[
(\polysub{B} f ) (a, k) \coloneqq (a, f \fcmp k).
\]
We say that an endofunctor on \catstdt is polynomial
if it is naturally isomorphic to a polynomial functor associated to a setoid family $B$.

In fact, the formalisation~\cite{Emmenegger2018a}
only contains the definition of the action of $\polysub{B}$
on setoids and functions.
\end{dfn}

We need to make sure that this definition does coincide
with the standard one  in~\eqref{eq:poly} for an arbitrary
locally cartesian closed category.

Let us begin by unfolding
the definition of $P_f$ in~\eqref{eq:poly}
for an extensional function $f \typing B \extfnc A$.
Recall that the three adjoint functors in~\eqref{eq:lccadj} are constructed
using $\sum$-types and $\prod$-types.
More precisely,
the first functor in~\eqref{eq:poly} maps a setoid $X$ to the projection
$\proj_2 \typing X \times B \extfnc B$,
and the last one maps a function $x \typing X \extfnc A$ to its domain $X$.
To describe the action of the functor $f_* \typing \catstdm/B \to \catstdm/A$ on
a function $y \typing Y \extfnc B$,
let us say that a pair $(a,g)$ in
\[
\sum_{a \typing A}
S(f)\fapp a \extfnc Y
\]
is a \emph{local section of $y$} if,
for every $(b,\_) \typing S(f)\fapp a$,
it is $y\fapp (g\fapp (b,\_)) \approx_B b$.
Two local sections $(a,g)$ and $(a',g')$ are equal if
there is $\alpha \typing a \approx_A a'$ such that
$g\fapp \approx g' \fcmp S(f)_{\alpha}$.
This defines a setoid $LS_y$ of local sections of $y$,
and the function $f_*\fapp y \typing LS_y \extfnc A$
is the obvious projection.

Let us now denote by $S(f)$ the setoid family constructed
from the extensional function $f$ in the proof of \cref{thm:stdfam}.
Unfolding the definitions, one sees that, for every $X$,
the setoid $P_f X$ is in bijection with $\polysub{S(f)} X$,
and that the bijection is natural in $X$.
If we denote by $E$ the inverse construction to $S$ from the same proof,
it also follows that for every setoid family $B \typing \sfam A$
the functor $P_{E(B)}$ is naturally isomorphic to $\polysub{B}$.

\begin{dfn}\label{def:alg}
Let $B \typing \sfam A$ be a setoid family.
\begin{enumerate}[beginpenalty=99,midpenalty=99]
\item\label{alg:map}
An \emph{algebra for $\polysub{B}$ in \catstdt} is
a setoid $X$ together with an extensional function $\polysub{B}X \extfnc X$.
\item\label{alg:mor}
Let $x \typing \polysub{B}X \extfnc X$ and $y \typing \polysub{B}Y \extfnc Y$
be algebras for $\polysub{B}$.
The \emph{setoid $\alg(x,y)$ of algebra morphisms from $x$ to $y$}
consists of functions $h \typing X \extfnc Y$
such that $h \circ x \approx y \circ \polysub{B}h$.
In formulas:
\[
|\alg(x,y)| \coloneqq \sum_{h \typing W \extfnc C}
h \fcmp x \approx y \fcmp (\polysub{B} h)
\]
and $(h,\_) \approx_{\alg(x,y)} (h',\_) \coloneqq h \approx h'$.
\end{enumerate}
\end{dfn}

\begin{dfn}\label{def:iscoh}
Let $B\typing \sfam A$ be a setoid family and $C$ a setoid.
A family of extensional functions
$F \typing \prod_{x \typing |A|} B\fapp x \extfnc C$
is \emph{coherent \wrt $B$} if, for every $a,a' \typing A$,
if $ \alpha\typing a \approx_A a'$ then $F\fapp a \approx F\fapp a' \circ B_{\alpha}$.

Let the type
\[
\isCoh{B}\fapp F \coloneqq
\prod_{(a,a'\typing A)} \prod_{(\alpha \typing a \approx a')}
	F\fapp a \approx F\fapp a' \circ B_{\alpha}
\]
be the type of proofs that $F$ is coherent.
\end{dfn}

\begin{lmm}\label{lem:extcoh}
Let $f \typing A' \extfnc A$ be an extensional function and
let $B \typing \sfam A$ be a setoid family over $A$.
For every family
$F \typing \prod_{x \typing |A'|} B\fapp (f\fapp x) \extfnc C$
of extensional functions, the following are equivalent.
\begin{enumerate}
\item
The family $F$ is coherent \wrt $f^*B \typing \sfam A'$.
\item
The function term
$\lambda x.(f\fapp x,F\fapp x) \typing A' \to \polysub{B}C$ is extensional.
\end{enumerate}
\end{lmm}

\begin{proof}
Let $\alpha' \typing a \approx_{A'} a'$.
Then $(f\fapp a,F\fapp a) \approx (f\fapp a',F\fapp a')$ \iff there is
$\alpha \typing f\fapp a \approx_A f\fapp a'$ such that
$F\fapp a \approx F\fapp a' \circ B_{\alpha}$.
Thus it clearly follows from coherence of $F$ taking
$\alpha \coloneqq f\extpf \alpha'$.
The converse is a direct consequence of proof irrelevance of $B$,
see~\eqref{eq:famirr} in \cref{rem:famirr}.
\end{proof}

\begin{rmk}
The value of the polynomial functor $\polysub{B}$ on a setoid $C$
can be described as the total setoid of the setoid family
$\lambda x. B\fapp x \extfnc C$ over $A$.
Accordingly, \cref{def:iscoh} and \cref{lem:extcoh}
can be phrased (and proved) in greater generality.
However we do not need that generality here.
\end{rmk}

\section{The algebra of extensional trees}
\label{sec:exttree}

Let $B \typing \sfam A$ be a setoid family over a setoid $A$.
In this section we construct an algebra from the W-type $\wty[|B|]$
on the underlying type family $|B| \typing |A| \to \ttunivm$.
We keep $A$ and $B$ fixed throughout the section and
we drop the subscript $B$ in the polynomial functor $\polysub{B}$.

The underlying setoid of the algebra is constructed in \cref{ssec:exttree}.
Before showing that $\wsup$ lifts to an algebra map in \cref{ssec:algmap},
we introduce the setoid family of immediate subtrees in \cref{ssec:treefam}.
\Cref{ssec:algmap} contains also a proof that the algebra map is invertible.

\subsection{The setoid of extensional trees}
\label{ssec:exttree}

To construct an algebra, we first need a setoid of trees,
that is, an equivalence relation on $\wty[|B|]$.
However, we should not expect $\wty[|B|]$
to be the underlying type of our setoid of trees though
or, in other words, we should not expect our relation to be total on $\wty[|B|]$.
Consider in analogy the situation described in \cref{rem:extper} with function terms:
the domain of the equivalence relation $\mathsf{extgen}$ on the function type consists precisely
of the extensional functions.

In order to find sufficient conditions for defining a suitable relation,
let us suppose that such a partial equivalence relation
$\approx_W \typing \wty[|B|] \to \wty[|B|] \to \ttunivm$
does exist.
This means that, if we denote by $W \typing \catstdm$
the setoid induced by $\approx_W$ on its domain
$|W| \coloneqq \sum_{w \typing \wty[|B|]} w \approx_W w$,
there is an initial algebra $ s \typing \poly{B} W \extfnc W$
whose underlying function term
$|s| \typing \left(\sum_{a \typing A} B\fapp a\extfnc W\right) \to |W|$
is a lift of
$\wsup \typing \left(\sum_{a \typing A}|B|\fapp a\to \wty[|B|]\right) \to \wty[|B|]$.
Since $s$ must be a bijection by Lambek's Theorem,
for every tree $w \jueq \wsup\fapp a\fapp f$ such that $w \approx_W w$ 
there is $g \typing B\fapp a \extfnc W$ such that $s(a,g) \approx_W w$.
This can be equivalently stated as follows.
\begin{equation}\tag{W1}\label{condW1}
\text{If $\wsup\fapp a\fapp f \approx_W \wsup\fapp a\fapp f$, then for every $b,b' \typing B\fapp a$,
if $b \approx_{B\fapp a} b'$ then $f\fapp b \approx_W f\fapp b'$.}
\end{equation}
Equivalently, for every tree $w$ in the domain of $\approx_W$
the immediate-subtree function term $\wb_w$ is extensional.
Second, again because $s$ is a bijection, and by definition of $\poly{B}$ in~\ref{def:poly},
it must be
\begin{equation}\tag{W2}\label{condW2}
\wsup \fapp a \fapp f \approx_W \wsup \fapp a' \fapp f'
\quad\text{\iff}\quad \begin{cases}
\,\text{there is $\alpha \typing a \approx_A a'$ such that}
\\
\,f \fapp b \approx_W f' \fapp (B_{\alpha}\fapp b)
\text{ for every $b \typing B\fapp a$.}
\end{cases}
\end{equation}

Let us now try to use (W1-2) to define $\approx_W$.
In condition \eqref{condW2} the relation $\approx_W$ 
occurs on the right-hand side only between immediate subtrees.
This feature makes \eqref{condW2} a possible candidate for the inductive step
in an inductive definition of the relation $\approx_W$.
However, it does not seem possible to derive \eqref{condW1} from \eqref{condW2}.
Instead, we combine them in the following definition.
Note that in \ref{def:Wper}.\ref{Wper:sub} below the equality between $b$ and $b'$ is over
the term $\alpha$ as in \eqref{def:catstd:depeq}.

\begin{dfn}\label{def:Wper}
Let $a,a' \typing A$, $f \typing |B|\fapp a \to \wty[|B|]$  and $f' \typing |B|\fapp a' \to \wty[|B|]$.
Two trees $w \jueq \wsup \fapp a \fapp f$  and $w' \jueq \wsup \fapp a' \fapp f'$  are \emph{matching} if
\begin{enumerate}
\item\label{Wper:node}
the nodes are equal, \ie
there is $\alpha \typing a \approx_A a'$, and
\item\label{Wper:sub}
for every $b \typing B\fapp a$  and $b' \typing B\fapp a'$,
if $b \approx_\alpha b'$ then the immediate subtrees
$f \fapp b$  and $f' \fapp b'$  are matching.
\end{enumerate}
A tree that is matching itself will be called
\emph{self-matching} or \emph{extensional}.

Accordingly, we define the type family
$\wper{B} \typing \wty[|B|] \to \wty[|B|] \to \ttunivm$
of proofs that two trees are matching as the (curried version of)
the dependent W-type on $\wty[|B|] \times \wty[|B|]$
with family of nodes
\[
N \coloneqq \lambda (w,w'). \wn \fapp w \approx_A \wn \fapp w'
\typing \wty[|B|] \times \wty[|B|] \to \ttunivm
\]
branching family
\[
Br \coloneqq \lambda (w,w'),\alpha. \sum_{b, b'} b \approx_\alpha b'
\typing \prod_I N\fapp (w,w') \to \ttunivm
\]
and arity function
\[
ar \coloneqq \lambda (w,w'),\alpha,(b,b',\beta). (\wb_w \fapp b, \wb_{w'} \fapp b')
\typing \prod_{I}\prod_{N} Br((w,w'),\alpha) \to I.
\]
\end{dfn}

The matching relation satisfies \eqref{condW1} and \eqref{condW2}:
for the former take $\alpha$ in~\ref{def:Wper}.\ref{Wper:node} to be reflexivity on $a$,
and for the latter take $b' \coloneqq B_{\alpha}\fapp b$ and reflexivity on $b'$
in~\ref{def:Wper}.\ref{Wper:sub}.

\begin{rmk}\label{rem:charWper}
For two trees $w,w' \typing \wty[|B|]$,
a proof term of \ref{def:Wper}.\ref{Wper:node} is
\[
\alpha \typing \wn \fapp w \approx_A \wn\fapp w'
\]
and a proof term of \ref{def:Wper}.\ref{Wper:sub} is
\[
\phi \typing \prod_{b \typing B\fapp (\wn \fapp w)}
\prod_{b' \typing B\fapp (\wn \fapp w')}
\Big(b \approx_\alpha b' \longrightarrow \wper{B} \fapp
(\wb_w \fapp b) \fapp (\wb_{w'} \fapp b')\Big).
\]

If we have two such terms, then a proof that $w$ and $w'$ are matching
is
\[
\wsupd\fapp (w,w')\fapp \alpha\fapp \phi \typing \wper{B}\fapp w\fapp w'.
\]

Conversely, if $w$ and $w'$ are matching with proof $\gamma$, then
$\alpha$ is obtained by applying the dependent node function
$\dwn_{(w,w')}$ from~\eqref{eq:dwn} to $\gamma$,
and $\phi$ by evaluating the dependent branching function
$\dwb_{(w,w')}$ from~\eqref{eq:dwb} on $\gamma$.
\end{rmk}

\begin{prp}\label{prop:wpersymtra}
The type family $\wper{B} \typing \wty[|B|]\to \wty[|B|] \to \ttunivm$
is a partial equivalence relation on $\wty[|B|]$,
that is, the following types are inhabited:
\begin{gather*}
\prod_{w, w' \typing \wty[|B|]} \wper{B} \fapp w \fapp w' \to \wper{B} \fapp w' \fapp w,
\\
\prod_{w, w', w'' \typing \wty[|B|]} \wper{B} \fapp w \fapp w' \to%
				    \wper{B} \fapp w' \fapp w'' \to \wper{B} \fapp w \fapp w''.
\end{gather*}
\end{prp}

\begin{proof}
The proof terms are obtained from straightforward applications
of the elimination rule for dependent W-types.
Alternatively, one may use \cref{rem:charWper} and recursion on $\wty[|B|]$.
\end{proof}

We can now form the setoid of extensional trees.

\begin{dfn}\label{def:exttrees}
Define the \emph{setoid $W_B$ of extensional trees} as the pair
$ (|W_B|, \approx_{W_B}) \typing \catstdm$,
where $|W_B|$ is the type of extensional (or self-matching) trees,
and two terms in $|W_B|$  are equal if their underlying trees are matching.
In formulas:
\[
|W_B| \coloneqq \sum_{w \typing \wty[|B|]} \wper{B} \fapp w \fapp w	\qquad \text{and} \qquad%
(w,\_) \approx_{W_B} (w',\_) \coloneqq \wper{B} \fapp w \fapp w'.\rule{0pt}{1.5em}
\]
We often leave the proof term in $\wper{B} \fapp w \fapp w$ implicit,
drop the subscript $B$ when it is clear from context,
and write $w \typing W$ to mean $w \typing \wty[|B|]$ and $w$ extensional.
\end{dfn}

Similarly to the situation described in \cref{rem:extper},
to prove that a tree $w$ is extensional it is sufficient to consider
an instance of the type of $\phi$ in \cref{rem:charWper},
namely the one where $\alpha$ is reflexivity at the node of $w$.

\begin{lmm}\label{lem:exttree}
A tree $w \jueq \wsup\fapp a \fapp f$ is extensional \iff
for every $b,b' \typing B\fapp a$, if $b \approx b'$ then
the trees $f\fapp b$ and $f\fapp b'$ are matching.

In particular, every immediate subtree of an extensional tree is extensional
and $f$ lifts to an extensional function $B\fapp a \extfnc W$.
\end{lmm}

\begin{proof}
If $w$ is extensional the conclusion follows by \eqref{condW1}.
The converse follows from \cref{rem:charWper} choosing reflexivity as $\alpha$.
\end{proof}

The function term $\wn \typing \wty[|B|] \to |A|$ lifts to
an extensional function $\wn \typing W \extfnc A$,
with proof of extensionality given by $\alpha$ in \cref{rem:charWper}.
By \cref{lem:exttree}, for every extensional tree $w$ the function term
$\wb_w \typing |B| \fapp (\wn \fapp w) \to \wty[|B|]$
lifts to an extensional function $\wb_w \typing B \fapp (\wn \fapp w) \extfnc W$.
So we have a family
\begin{equation}\label{eq:treeonbr}
\wb \typing \prod_{w \typing W} B \fapp (\wn \fapp w) \extfnc W
\end{equation}
of extensional functions.
We continue writing $\wb_w$ for $\wb\fapp w$ and
we will do the same also
for others dependent function terms on $W$.
As the subscript will always be a tree in $W$,
no confusion should arise with the action of setoid families on equalities.

Also to prove that two \emph{extensional} trees are matching
it is enough to consider
just an instance of the type of $\phi$ in \cref{rem:charWper},
as for extensional functions in \cref{rem:extper}.

\begin{lmm}\label{lem:charWeq}
Two extensional trees
$ w \jueq \wsup\fapp a\fapp f$  and $w' \jueq \wsup\fapp a'\fapp f'$  are matching \iff
there is $\alpha \typing a \approx_A a'$ such that $f \approx f' \circ B_{\alpha}$.

Equivalently, for every $w, w' \typing W$,
it is $w \approx_W w'$ \iff there are
\[
\alpha \typing \wn \fapp w \approx_A \wn\fapp w'
\hspace{3em}\text{and}\hspace{3em}
\psi \typing \wb_w \approx \wb_{w'} \circ B_{\alpha}.
\]
\end{lmm}

\begin{proof}
If $w$ and $w'$ are matching, the conclusion follows by \eqref{condW2}.
Conversely, by \cref{rem:charWper} it is enough
to construct a term
\[
\phi \typing \prod_{b,b'}
b \approx_{\alpha} b' \to f\fapp b \approx_W f'\fapp b'.
\]
By \cref{lem:exttree} we may assume
that the function $f'$ is extensional.
It follows that
for every $b \typing B\fapp a$, $b' \typing B\fapp a'$,
if $\beta \typing b \approx_{\alpha} b'$,
then $f\fapp b \approx_W f'\fapp b'$
using
$\psi\fapp b \typing f\fapp b \approx_W f'\fapp (B_{\alpha}\fapp b)$
and $f'\extpf \beta$.
\end{proof}

\subsection{Setoid families on extensional trees}
\label{ssec:treefam}

Extensionality of $\wn \typing W \extfnc A$ yields,
for every $\gamma \typing w \approx_W w'$ between extensional tree,
a function
\begin{equation}\label{eq:brtrsp}
B_{\wn\extpf \gamma} \typing B\fapp (\wn\fapp w) \extfnc B\fapp (\wn\fapp w')
\end{equation}
where $\wn\extpf \gamma$ denotes the proof of
extensionality~\eqref{eq:extfun} of $\wn$ applied to $\gamma$.
Since $B$ is a setoid family,
the family of setoids
$\lambda w. B\fapp (\wn\fapp w) \typing W \to \catstdm$
is a setoid family with the functions in~\eqref{eq:brtrsp}
as transports.
We refer to this setoid family as
the family of \emph{branches of trees}.

There is also another natural setoid family on extensional trees
which will turn out to be instrumental in the characterisation of
algebra morphisms in \cref{sec:init}.

\begin{dfn}\label{def:imsub}
Let $w \jueq \wsup\fapp a\fapp f$ be an extensional tree.
An \emph{immediate subtree of $w$} is a subtree index
$b \typing B\fapp a$,
and two immediate subtrees $b$  and $b'$ of $w$  are equal if the
corresponding trees $f\fapp b$  and $f\fapp b'$  are matching.

Accordingly, the setoid of \emph{immediate subtrees of $w$},
denoted $\isbtree \fapp w$,
has $|B| \fapp (\wn \fapp w)$  as underlying type,
and 
\[
s \approx_{\isbtree \fapp w} s' \coloneqq \wb_w \fapp s \approx_W \wb_w \fapp s'
\]
as equality.

For $s \typing \isbtree \fapp w$  and $\gamma \typing w \approx_W w'$,
the assignment
\[
\isbtree_\gamma \fapp s \coloneqq B_{\wn\extpf \gamma} \fapp s
\typing \isbtree\fapp w'
\]
defines transport maps for $\isbtree$.
Hence we obtain a setoid family $\isbtree \typing \sfam W$,
the \emph{family of immediate subtrees}.
\end{dfn}

\begin{rmk}\label{rem:imgfact}
The category of setoids \catstdt is exact and,
assuming the identity type,
it is the exact completion of its subcategory of discrete setoids.
We refer to~\cite{EmmeneggerPalmgren2017} for details
but see also~\cite{vdBergMoerdijk2018}.
Being exact, every function has a factorisation
as a regular epi followed by a mono, called image factorisation.
Being an exact completion, every function $f \typing X \extfnc Y$
has a canonical such factorisation.
Briefly: the setoid $Z$ has $|X|$  as underlying type and
$\lambda x,x'.\left(f\fapp x\approx_Y f\fapp x'\right)$  as equality;
the regular epi $X \extfnc Z$ has $\lambda x.x$ as underlying function term,
which is extensional because $f$ is;
and the mono $Z\extfnc Y$ has $|f|$ as underlying function term,
which is extensional by definition of $Z$.
In fact, this canonical factorisation
exists also when identity types are not assumed.

By applying this factorisation to the function
$ \wb_w \typing B\fapp (\wn\fapp w) \extfnc W$ 
we recover the setoid of immediate subtrees of $w$:
\[\xycenterm{
B \fapp (\wn \fapp w)	\ar@{=>}[dr]_{\istq_w} \ar@{=>}[rr]^{\wb_w}
&&	W
\\
&	\isbtree \fapp w	\ar@{=>}[ur]_{\swb_w}	&
}\]
where $\istq_w \coloneqq (\mathsf{id},\wb_w\extpf)$
and $\swb_w \coloneqq (|\wb_w|,\mathsf{id})$
denote the regular epi and the mono, respectively,
arising from the factorisation of $\wb_w$.

These form in turn two families of extensional functions
\begin{equation}\label{eq:estree}
\istq \typing \prod_{w \typing W} B\fapp (\wn\fapp w) \extfnc \isbtree\fapp w
\qquad\text{and}\qquad
\swb \typing \prod_{w \typing W} \isbtree\fapp w \extfnc W.
\end{equation}
The family $\istq$ is coherent in the sense that
\begin{equation}\label{eq:istqcoh}
\isbtree_{\gamma} \circ \istq_w \approx \istq_{w'} \circ B_{\wn\extpf \gamma}
\end{equation}
holds since
$\isbtree_{\gamma}$ and $B_{\wn\extpf \gamma}$ have the same function term.
The family $\swb$ is coherent in the sense of \cref{def:iscoh} by \cref{rem:cohbr}.
\end{rmk}

\begin{rmk}\label{rem:cohbr}
By \eqref{condW2},
the family of functions $\wb$
from~\eqref{eq:treeonbr} is coherent \wrt the setoid family $\isbtree$
in the sense of \cref{def:iscoh},
\ie for every $\gamma \typing w \approx_W w'$ there is
\[
\wbcoh \fapp \gamma \typing
\wb_w \approx \wb_{w'} \circ \isbtree_{\gamma}.
\]
The proof term $\wbcoh \fapp \gamma$
is obtained instantiating $\phi$ in \cref{rem:charWper}
with $b' \coloneqq \isbtree_{\gamma}\fapp b$  and reflexivity on $b'$.

Similarly, the family of functions $\swb$
from~\eqref{eq:estree} is coherent \wrt the setoid family $\isbtree$
in the sense of \cref{def:iscoh},
\ie for every $\gamma \typing w \approx_W w'$ there is
\begin{equation}\label{eq:streecoh}
\swbcoh\fapp \gamma \typing \swb_w \approx \swb_{w'} \circ \isbtree_{\gamma}.
\end{equation}
The proof term $\swbcoh\fapp \gamma$ is obtained from $\wbcoh$
using~\eqref{eq:istqcoh} and the fact that $\istq_w$ is a (canonical regular) epi.
\end{rmk}

\subsection{Construction of the algebra map}
\label{ssec:algmap}

We construct the algebra map by showing that
the function term $\wsup$  lifts to an extensional function
$\ws \typing \poly{B} W \extfnc W$.

For $a \typing A$  and $f \typing B \fapp a \extfnc W$,
let $f_0 \coloneqq \proj_1 f \typing |B| \fapp a \to \wty[|B|]$.
The tree $w \coloneqq \wsup\fapp a\fapp f_0$  is extensional by \cref{lem:exttree},
so there is a function term
\begin{equation}\label{eq:amtrm}
\ws \typing \poly{B}W \to W.
\end{equation}
Given explicitly by
\begin{equation}\label{eq:ameq}
\ws \fapp (a,f) \coloneqq (\wsup \fapp a \fapp f_0,%
    \wsupd \fapp (w,w) \fapp (\rho \fapp a) \fapp \epsi'_f) \typing W
\end{equation}
where $\rho$  is reflexivity of $A$, $\epsi_f$  is obtained uncurrying
$ f\extpf \typing \prod_{b,b'} b \approx_{B \fapp a} b' \to f\fapp b \approx_W f\fapp b'$ 
and $\wsupd \fapp (w,w) \fapp (\rho \fapp a) \fapp \epsi_f \typing \wper{B} \fapp w \fapp w$ 
is a proof that $w$  is extensional.

\begin{lmm}
The function term $\ws \typing \poly{B} W \to W$ is extensional.
\end{lmm}

\begin{proof}
By \cref{def:poly}, to have two equal elements in the domain is to have
$a, a' \typing A$, $f \typing B \fapp a \extfnc W$, $f' \typing B \fapp a' \extfnc W$,
and a term $\alpha \typing a \approx a' $such that $f \approx f' \fcmp B_\alpha$.
\Cref{lem:charWeq} and \cref{def:Wper} yield the claim.
\end{proof}

We conclude this section with a proof that the algebra map
$\ws \typing \poly{B}W \extfnc W$ is invertible.

\begin{prp}\label{prop:siso}
There is $\wus \typing W \extfnc \poly{B} W$ such that
$\ws \fcmp \wus \approx id_W$  and $\wus \fcmp \ws \approx id_{\poly{B} W}$.
\end{prp}

\begin{proof}
The function $\wus$ maps an extensional tree $w \jueq \wsup\fapp a\fapp f_0$ to
the pair $(a,f) \typing \poly{B}W$, where $f \typing B\fapp a \extfnc W$
is the extensional function of $w$ from \cref{lem:exttree}.
In formula
\[
\wus \fapp w \coloneqq (\wn \fapp w, \wb_w) \typing \poly{B}W.
\]
It is extensional because $\wn$ is extensional,
and because $\wb$ is coherent by~\cref{rem:cohbr}.
The equation $\wus \circ \ws \approx \mathsf{id}_{\poly{B}W}$
is straightforward from the definitions.
For $\ws \circ \wus \approx \mathsf{id}_W$
use \cref{lem:charWeq}.
\end{proof}

\section{Initiality of the algebra of extensional trees}
\label{sec:init}

This section is mainly devoted to the characterisation of algebra morphisms
as telescopic functions in \cref{thm:commchar}.
The existence of W-types in \catstdt follows in \cref{thm:init}.

\subsection{Characterisation of algebra morphisms}
\label{ssec:charalgmor}

Throughout this section we fix 
setoids $A$ and $C$, a setoid family $B \typing \sfam A$
and a $\polysub{B}$-algebra $a_C \typing \polysub{B} C \extfnc C$.
As in the previous section, we drop the index $B$ from
the polynomial functor $\polysub{B}$ and the algebra of extensional trees $W_B$.

We begin with a simple observation.
Recall from \cref{def:alg}.\ref{alg:mor}
that $\alg(\ws,a_C)$ denotes the setoid of algebra morphisms from
$\ws \typing \poly{B}W \extfnc W$ to $a_C \typing \poly{B}C \extfnc C$,
\ie those functions $h \typing W \extfnc C$ such that
\begin{equation}\label{eq:algcomm}
h \circ \ws \approx a_C \circ \poly{B}h.
\end{equation}

Every function $h \typing W \extfnc C$ gives rise
to a family
\begin{equation}\label{eq:gluedata}
F \typing \prod_{w \typing W} \isbtree\fapp w \extfnc C
\end{equation}
defined by $F_w \coloneqq h \circ \swb_w$.
It follows from~\eqref{eq:algcomm}
that the function $h$ is an algebra morphism \iff,
for every $w \jueq \ws(a,f)$,
\begin{equation}\label{eq:glueptw}
h\fapp w
\approx a_C(a, F_w \circ \istq_w).
\end{equation}
In particular, the value of the algebra morphism $h$
at the tree $w$ is obtained by ``glueing''
the values of $h$ on the immediate subtrees of $w$,
collected in $F_w$,
according to the algebra $a_C$.

Note that the right-hand side of~\eqref{eq:glueptw} makes sense
for any family $F$ as in~\eqref{eq:gluedata}.
Let us say that a function $k \typing \isbtree\fapp w \extfnc C$
is a \emph{glueing datum in $C$ for $w$},
and that a family $F$ as in \eqref{eq:gluedata}
is a \emph{family \famgldtin{C}}.
Two families $F$ and $F'$ \famgldt are equal if,
for every $w \typing W$, the glueing data for $w$ are equal,
\ie $F_w \approx F'_w$.
Given a family $F$ \famgldtin{C} as in~\eqref{eq:gluedata},
we can use~\eqref{eq:glueptw} to define a function term $h \typing W \to C$.
\Cref{prop:restrglue} isolates conditions on $F$
for this term to be extensional and,
further, an algebra morphism.

The problem of constructing an element of $\alg(\ws,a_C)$
could thus be reduced to the problem of finding a suitable family
\famgldtin{C} and to prove that it is the unique such.
We do so in \cref{thm:commchar},
where we prove that $\alg(\ws,a_C)$ is in bijection with
a subsetoid of families \famgldtin{C} 
that we name, in \cref{def:tel}, \famstel.
In \cref{prop:telglue} we construct
a \famtel by direct induction on the underlying W-type
$\wty[|B|]$ of $W$.
Its uniqueness is in \cref{cor:tlabssing}.
It follows that $\alg(\ws,a_C)$ is a singleton.
This is in accordance with the result in~\cite{AGS2017},
that an algebra $a$ is equivalent to a W-type
if and only if the type of morphisms of algebras out of $a$ is contractible,
once we observe that, in the setoid interpretation,
a type is contractible precisely
when the corresponding setoid is a singleton.

There is no reason, at this stage, to prefer families of functions
on immediate subtrees as in~\eqref{eq:gluedata} over
families of functions $B\fapp (\wn\fapp w) \extfnc C$
on branches of trees,
but our strategy does not seem to work with the latter.
One motivation is discussed in \cref{rem:whysub}.

First, we give names to the two constructions outlined above.

\begin{dfn}\label{def:restrglue}
Let $h \typing W \extfnc C$ be an extensional function.
The \emph{restriction of $h$ to (immediate) subtrees}
is the family \famgldt
\begin{equation}\label{eq:restr}
\rest\fapp h \coloneqq \lambda w. h\circ \swb_w
	\typing \prod_{w \typing W} \isbtree\fapp w \extfnc C.
\end{equation}

Let $F \typing \prod_{w:W} \isbtree\fapp w \extfnc C$ be a family \famgldt.
The \emph{glueing of $F$ (along $a_C$)} is the function term
\begin{equation}\label{eq:glue}
\glu\fapp F \coloneqq \lambda w. a_C(\wn\fapp w, F_w \circ \istq_w) \typing W \to C.
\end{equation}
\end{dfn}

We need one more definition to state the next result.
Coherence for a family of extensional functions was defined in~\ref{def:iscoh}.

\begin{dfn}\label{def:lab}
Let us denote by $\submaps$ the setoid of families \famgldtin{C}.
The underlying type is $\prod_{w \typing W}\isbtree\fapp w \extfnc C$
and the equality between two families $F$ and $F'$ is
\[
F \approx F' \ \coloneqq\
\prod_{w \typing W} F_w \approx_{(\isbtree\fapp w \extfnc C)} F'_w.
\]

The subsetoid of $\submaps$ on the coherent families \famgldt is
denoted $\cLab$ and has underlying type
\[
|\cLab| \coloneqq \sum_{F \typing \submaps} \isCoh{\isbtree}\fapp F.
\]
\end{dfn}

\begin{prp}\label{prop:restrglue}\hfill
\begin{enumerate}
\item\label{restrglue:restrcoh}
For every function $h \typing W \extfnc C$,
the family \famgldt
$\rest\fapp h$ defined in \eqref{eq:restr} is coherent.
\item\label{restrglue:cohgluext}
For every coherent family \famgldt $F$,
the function $\glu\fapp F$ defined in \eqref{eq:glue} is extensional.
\item\label{restrglue:restrext}
The function $\rest \typing (W \extfnc C) \to \cLab$ is extensional.
\item\label{restrglue:gluext}
The function $\glu \typing \cLab \to (W \extfnc C)$ is extensional.
\item\label{restrglue:grid}
For every $h \typing W \extfnc C$
\[
\glu\fapp (\rest\fapp h) \approx h
\,\longleftrightarrow\,
h \circ \ws \approx a_C \circ \poly{B}h.
\hspace{10ex}
\]
\item\label{restrglue:rgid}
For every coherent family
$F \typing \prod_{w:W} \isbtree\fapp w \extfnc C$ \famgldt
\[
\rest\fapp (\glu\fapp F) \approx F
\,\longrightarrow\,
\glu\fapp F \circ \ws \approx a_C \circ \poly{B}(\glu\fapp F).
\]
\end{enumerate}
\end{prp}

\begin{proof}
\eqref{restrglue:restrcoh}
For every $\gamma \typing w \approx_W w'$, it is
$h \circ \swb_w \approx h \circ \swb_{w'} \circ \isbtree_{\gamma}$
since $\swb_w$ is coherent by \eqref{eq:streecoh} in \cref{rem:cohbr},
and $h$ is extensional.

\eqref{restrglue:cohgluext}
This is \cref{lem:extcoh} together with extensionality of $a_C$.

\eqref{restrglue:restrext}
If $h \approx h'$, then for every $w \typing W$ it is
$h \circ \swb_w \approx h' \circ \swb_w$
by extensionality of $h$.

\eqref{restrglue:gluext}
Let $\varphi \typing F \approx F'$.
If $F_w \approx F'_w$, then by extensionality of $a_C$
and proof irrelevance of $B$ from~\eqref{eq:famirrref} it is
\[
a_C(a, F_w) \approx_C a_C(a, F'_w).
\]

\eqref{restrglue:grid}
The claim follows from the fact that,
for every $w \jueq \ws(a,f)$, it is
\begin{equation}\label{eq:amaltcomm}
\begin{split}
\glu\fapp (\rest\fapp h)\fapp w &\jueq a_C(a,h\circ \swb_w \circ \istq_w)
\\&\approx
a_C(a,h \circ f)
\\&\approx
a_C \circ \poly{B}h (a,f)
\\&\approx
a_C \circ \poly{B}h \circ \wus\fapp w
\end{split}
\end{equation}
by, in order, \eqref{eq:restr} and \eqref{eq:glue}, \cref{rem:imgfact},
definition of polynomial functor in~\ref{def:poly}
and definition of $\wus$ in \cref{prop:siso}.

\eqref{restrglue:rgid}
Unfolding the right hand side for $w \jueq \ws(a,f)$ yields
\[
a_C(a, F_w \circ \istq_w) \approx_C
a_C(a,\lambda b.a_C(\wn\fapp (f\fapp b), F_{(f\fapp b)}\circ \istq_{(f\fapp b)}))
\]
By extensionality of $a_C$ and proof irrelevance~\eqref{eq:famirrref},
this holds if for every $b \typing B\fapp a$ it is
\[\begin{split}
F_w\fapp b \approx_C a_C(\wn\fapp (f\fapp b), F_{(f\fapp b)}\circ \istq_{(f\fapp b)}).
\end{split}\]
Now, by definitions~\eqref{eq:restr} and~\eqref{eq:glue} it is
\[
\rest\fapp (\glu\fapp F)\fapp w\fapp b
\jueq a_C(\wn\fapp (f\fapp b),F_{(f\fapp b)}\circ \istq_{(f\fapp b)})
\]
and the claim follows.
\end{proof}

\begin{crl}\label{cor:restrglue}
The functions $\rest$ and $\glu$ establish a bijection between
the setoid $\alg(\ws,a_C)$ and those coherent families $F$ \famgldtin{C}
such that $\rest\fapp (\glu\fapp F) \approx F$.
\end{crl}

\begin{proof}
\Cref{prop:restrglue}.\ref{restrglue:grid}
entails that $\rest\fapp (\glu\fapp (\rest\fapp h)) \approx \rest\fapp h$
for every algebra morphism $h$.
Conversely, If $F$ is such that $\rest\fapp (\glu\fapp F) \approx F$,
then $\glu\fapp F$ is in $\alg(\ws,a_C)$
by \cref{prop:restrglue}.\ref{restrglue:rgid}.
\end{proof}

The bijection in \cref{cor:restrglue} suggests
a recursive characterisation of algebra morphisms into $C$
as those morphisms obtained by glueing a coherent family \famgldtin{C},
such that each component $F_w$ is itself obtained by glueing
a suitable family, and so on.
In order to make this intuition precise,
we need to bring coherent families \famgldt
and glueing one ``subtree level'' up.
We do not need to go any higher thanks to dependent W-types,
which will cover the other cases for us in \cref{def:tel},
see also \cref{rem:telalt}.
The characterisation is accomplished in \cref{thm:commchar}.

For $w$ an extensional tree, we say that a family
\[
G \typing \prod_{s \typing \isbtree \fapp w}
\isbtree \fapp (\swb_w \fapp s) \extfnc C.
\]
is a \emph{family \famglsubdt for $w$ in $C$}.
We write $G_s$ for $G\fapp s$ as we do for families \famgldt on subtrees,
and we refer to terms of type $\isbtree \fapp (\swb_w \fapp s)$
as \emph{2-subtrees of $w$}.

By extensionality of $\swb_w$, the family of setoids
$\lambda b. \isbtree\fapp (\swb_w\fapp b)$ over $\isbtree\fapp w$
is a setoid family in the same way
as for the family of subtrees $\isbtree$ in \cref{def:imsub}.
We can thus form the setoid $\csLab\fapp w$ of
\emph{coherent families \famglsubdt for $w$},
as we did in \cref{def:lab} for coherent families \famgldt of subtrees.

In fact, $\csLab$ is another setoid family over $W$,
whose transport function
\[
\csLab_\gamma \typing \csLab \fapp w \to \csLab \fapp w'
\]
for $\gamma \typing w \approx_W w' $
is defined on $G$ and $s \typing \isbtree \fapp w'$ as
\begin{equation}\label{eq:cohssmaptrsp}
\csLab_\gamma\fapp G \fapp s \coloneqq G_{(\isbtree_{\gamma^{-1}} \fapp s)}
\fcmp \isbtree_{(\swbcoh \fapp \gamma^{-1}\fapp s')}
\end{equation}
where $\swbcoh \fapp \gamma^{-1}\fapp s' \typing \swb\fapp s' \approx \swb\fapp (\isbtree_{\gamma^{-1}} \fapp s')$.
It follows that
\begin{equation}\label{eq:cohssmapeq}
G \approx_\gamma G'
\ \longleftrightarrow\
\prod_{s \typing \isbtree \fapp w}\hspace{-.5em}
G_s	\approx
      G'_{(\isbtree_\gamma \fapp s)} \fcmp \isbtree_{(\swbcoh\fapp \gamma\fapp s)}
\end{equation}
for all
$G \typing \csLab\fapp w$ and
$G' \typing \csLab\fapp w'$.

As in~\eqref{eq:glue} for every family $G$ \famglsubdt for $w$
there is a function term
\[
\sglu_w\fapp G \typing \isbtree\fapp w \to C
\]
defined on $s \typing \isbtree\fapp w$ by
\begin{equation}\label{eq:subglue}
\sglu_w\fapp G \fapp s \coloneqq a_C \fapp \left(\wn \fapp (\swb_w\fapp s), G_s \fcmp \istq_{(\swb_w\fapp s)}\right).
\end{equation}

\begin{lmm}\label{lem:subglue}\hfill
\begin{enumerate}
\item\label{subglue:valext}
For every coherent family \famglsubdt $G \typing \csLab\fapp w$,
the function $\sglu_w\fapp G$ defined in~\eqref{eq:subglue} is extensional.
\item\label{subglue:ext}
For every $w \typing W$,
the function $\sglu_w \typing \csLab\fapp w \to (\isbtree\fapp w \extfnc C)$
is extensional.
\item\label{subglue:coh}
For every $\gamma \typing w  \approx_W w'$, $G \typing \csLab \fapp w $
and $G' \typing \csLab \fapp w' $
there is
\[
\sglucoh \fapp \gamma \typing G \approx_\gamma G'	\longrightarrow%
	  \sglu_w\fapp G \approx \sglu_{w'}\fapp G' \fcmp \isbtree_{\gamma}.
\]
\end{enumerate}
\end{lmm}
\begin{proof}
\eqref{subglue:valext}
It follows from \cref{lem:extcoh} and extensionality of $a_C$.

\eqref{subglue:ext}
It follows from~\eqref{subglue:coh} with reflexivity as $\gamma$.

\eqref{subglue:coh}
Let $\gamma$, $G$ and $G'$ be as above,
and let $\varphi \typing G \approx_\gamma G'$ and $s \typing \isbtree \fapp w$.
Using extensionality of $a_C$, it is enough to show
\[
\left(\wn \fapp (\swb_w\fapp s), G_s \right) \approx%
\left(\wn \fapp (\swb_{w'}\fapp (\isbtree_\gamma \fapp s)),%
  G'_{(\isbtree_\gamma \fapp s)}\right).
\]
For the first component we may take
$\alpha \coloneqq \wn\extpf (\swbcoh\fapp \gamma)$
from \cref{rem:cohbr}.
It remains to show that, for every $s \typing \isbtree\fapp w$,
\[
G_s	\approx
      G'_{(\isbtree_\gamma \fapp s)} \fcmp \isbtree_{(\swbcoh \fapp \gamma)}.
\]
By~\eqref{eq:cohssmapeq} it is enough to use $\varphi$.
\end{proof}

\begin{rmk}\label{rem:whysub}
\Cref{lem:subglue} is the reason for choosing to work with
the setoid family of immediate subtrees instead of the setoid family of branches of trees
from \cref{ssec:treefam}.

Indeed, \cref{prop:restrglue} works equally well with immediate subtrees
replaced by branches of trees, in the sense that there are function terms
between extensional functions $W \extfnc C$ and families \famgldt on branches
satisfying the same properties, mutatis mutandis, of $\rest$ and $\glu$.
Problems arise at the next subtree level.

If we say that, for an extensional tree $w$, a term of type
\[
\prod_{b \typing B\fapp (\wn\fapp w)}
B\fapp (\wn\fapp (\wb_w\fapp b)) \extfnc C
\]
is a \emph{family \famglsubdt on branches for $w$ in $C$},
we can then form, for every $w \typing W$,
the setoid $CB_w$ of coherent families \famglsubdt on branches for $w$.
Nevertheless,
it seems that transport of families \famglsubdt on branches over
$\gamma \typing w\approx_W w'$ does not preserve coherence.
It follows that, contrary to $\csLab$,
the family $\lambda w.CB_w$ is not a setoid family on $W$.
In particular, there is no function term like
$\csLab_\gamma$ from~\eqref{eq:cohssmaptrsp},
which is needed, together with \cref{lem:subglue}.\ref{subglue:coh},
in a crucial step in the proof of \cref{lem:teltransp}.
The issue seems to arise because branches, contrary to subtrees,
do not form a subsetoid of~$W$.
\end{rmk}

Note that there is no version of $\rest$ for 2-subtrees:
a function $k \typing \isbtree\fapp w \extfnc C$ cannot be restricted
to $\isbtree\fapp (\swb_w\fapp s)$ for some $s \typing \isbtree\fapp w$
simply because an immediate subtree of $\swb_w\fapp s$
is not an immediate subtree of $w$, that is,
being an immediate subtree is not a transitive property.
This entails that, in order to say that $k$ is a glueing,
we need to assume that a coherent family \famglsubdt for $w$ exists.

\begin{dfn}\label{def:tel}
Let $w \jueq \ws(a,f)$ be an extensional tree.
A function $k \typing \isbtree\fapp w \extfnc C$ is
\emph{telescopic over $w$} if
\begin{enumerate}
\item\label{tel:lab}
the function $k$ is the glueing of
a coherent family \famglsubdt $G \typing \csLab\fapp w$,
\ie $k \approx \sglu_w\fapp G$, and
\item\label{tel:tel}
for every $s \typing \isbtree\fapp w$,
the function $G_s \typing \isbtree\fapp (f\fapp s) \extfnc C$ is telescopic over $f\fapp s$.
\end{enumerate}
Accordingly, we define the type family
\[
\isTel \typing \left(\sum_{w \typing W}\isbtree\fapp w \extfnc C\right) \to \ttunivm
\]
of proofs that a function on immediate subtrees is telescopic
as the dependent W-type on
$I \coloneqq \sum_{w \typing W} (\isbtree \fapp w \extfnc C)$
with family of nodes
\[
N \coloneqq \lambda (w,k). \sum_{G \typing \csLab \fapp w}
	k \approx \sglu_w\fapp G
\typing \left(\sum_{w \typing W} \isbtree \fapp w \extfnc C\right) \to \ttunivm
\]
branching family
\[
Br \coloneqq \lambda (w,k),(G,\_).\isbtree \fapp w
\typing \prod_{(w,k) \typing I} N\fapp (w,k) \to \ttunivm
\]
and arity function
\[
ar \coloneqq \lambda (w,k),(G,\_),s.(\swb_w\fapp s, G_s)
\typing \prod_{I}\prod_{N} \isbtree\fapp w \to
\left(\sum_{w \typing W} \isbtree \fapp w \extfnc C\right).
\]

Henceforth we write $\isTel_w\fapp k$ for $\isTel\fapp (w,k)$
and we may drop the subscript $w$ if it is clear from context.
\end{dfn}

\begin{dfn}\label{def:telstd}
A family $F$ \famgldtin{C} is \emph{telescopic} if
each of its components is a telescopic function over $w$, that is, if
\[
\prod_{w \typing W} \isTel\fapp F_w
\]
is inhabited.
We say `telescopic family' to mean
`telescopic family \famgldtin{C}'.

Telescopic families \famgldt form a subsetoid $\tLab$ of $\submaps$.
\end{dfn}

\begin{rmk}\label{rem:telalt}
If we say, similarly,
that a family $G$ \famglsubdt for $w$ is telescopic over $w$
if each component $G_s$ is telescopic over
the immediate subtree $\swb_w\fapp s$,
then \cref{def:tel} can be phrased as:
a function is telescopic over $w$
if it is the glueing of a telescopic family \famglsubdt for $w$.
\end{rmk}

The canonical components of a proof that a function
$k \typing \isbtree\fapp w \extfnc C$
is telescopic are as follows,
see also \cref{dwrules}.
From the dependent node function~\eqref{eq:dwn} we obtain,
for every $T \typing \isTel_w \fapp k$, a coherent family \famglsubdt
\begin{equation}\label{eq:cslab}
\slab \fapp k\fapp T \typing \csLab \fapp w
\end{equation}
and a proof
\begin{equation}\label{eq:sglue}
\telsglu \fapp k\fapp T \typing k \approx \sglu_w\fapp (\slab \fapp k \fapp T),
\end{equation}
that $k$ is a glueing of the family $\slab\fapp k\fapp T$.
The dependent subtree function~\eqref{eq:dwb} yields a proof
\begin{equation}\label{eq:cslabTel}
\telslab \fapp k\fapp T \typing
  \prod_{(s \typing \isbtree \fapp w)}%
    \isTel_{(\swb_w\fapp s)} \fapp (\slab \fapp k \fapp T \fapp s).
\end{equation}
that every function in the family $\slab\fapp k\fapp T$ is telescopic,
that is, the family $\slab$ \famglsubdt is itself telescopic,
see~\cref{rem:telalt}.

Clearly, terms in~\eqref{eq:cslab} and~\eqref{eq:sglue} are proof terms for
\ref{def:tel}.\ref{tel:lab}, and the term in~\eqref{eq:cslabTel}
is a proof term for \ref{def:tel}.\ref{tel:tel}.

As we will extensively use the elimination principle in \cref{dwrules}
of the dependent W-type of telescopic functions from~\cref{def:tel},
let us unfold its inductive hypothesis too.
Let $V \typing \prod_{w,k} \isTel\fapp k \to \ttunivm$
be a (small) type family.
The inductive hypothesis tells us that, for
$T \jueq \wsupd\fapp (w,k)\fapp G\fapp E \typing \isTel\fapp k$
and for every $s \typing \isbtree\fapp w$,
there is a term
\begin{equation}\label{eq:telIH}
I\!H\fapp s \typing V\fapp (\swb_w\fapp s)\fapp G_s\fapp (E\fapp s)
\end{equation}
where $G \jueq \slab\fapp T$ is the coherent family
\famglsubdt obtained from $T$ as in~\eqref{eq:cslab},
and the term $E \jueq \telslab\fapp T\fapp$ is the proof~\eqref{eq:cslabTel}
that every $G_s$ is telescopic over $\swb_w\fapp s$.

The fact that functions $\isbtree\fapp w \extfnc C$
cannot be restricted to families \famglsubdt,
discussed right before \cref{def:tel},
does not prevent us from showing that
restrictions of algebra morphisms are telescopic,
since a function $W \extfnc C$ can be restricted to any ``subtree level''.
In particular, to immediate subtrees and 2-subtrees.

\begin{lmm}\label{lem:restristel}
Let $h \typing W \extfnc C$ be an algebra morphism.
Then for every $w \typing W$, the function
\[
\rest\fapp h\fapp w \coloneqq h \fcmp \swb_w \typing \isbtree\fapp w \extfnc C
\]
defined in~\eqref{eq:restr} is telescopic over $w$.
\end{lmm}

\begin{proof}
The proof is by $\wty$-elimination on $w \typing \wty[|B|]$ into the type
\[
\wper{B} \fapp w \fapp w	\,\ \longrightarrow%
\isTel \fapp (h \fcmp \swb_w).
\]
Note that, whenever we wish to
apply the inductive hypothesis to a subtree index of $w$,
we have to provide a proof that the subtree is extensional.
Such a proof will always be available since
$w$ is extensional by assumption and
immediate subtrees of extensional trees are extensional by \cref{lem:exttree}.
So we may assume without loss of generality
that all the trees that we will be dealing with are extensional.

Let then $w \jueq \wsup \fapp a \fapp f$.
By \cref{def:tel}, in order to apply $\wsupd$,
we first need to provide a coherent family \famglsubdt $G \typing \csLab \fapp w$ and
a proof that $h \fcmp \swb_w \approx \sglu_w\fapp G$.

The family $G$ \famglsubdt consists of the restrictions of $h$
to 2-subtrees of $w$,
namely
\[
G_s \coloneqq h \fcmp \swb_{(\swb \fapp s)}
\typing \isbtree\fapp (\swb_w\fapp s) \extfnc C
\]
for every $s \typing \isbtree\fapp w$.
Its coherence follows as in~\ref{prop:restrglue}.\ref{restrglue:restrcoh}
from coherence of $\swb$  and extensionality of $h$.
For every $s \typing \isbtree \fapp w$ it is
\[\begin{split}
h\fapp (\swb_w \fapp s)
&\approx	a_C \fcmp (\poly{B} h)\fcmp \wus \fapp (\swb_w\fapp s)
\\&\approx
a_C \fapp (\wn \fapp (\swb_w\fapp s), h \fcmp \swb_{(\swb \fapp s)})
\\&\approx
\sglu_w\fapp G \fapp s
\end{split}\]
by the fact that $h$ is an algebra morphism, definition of $\wus$ and $\poly{B}$,
and definition of $\sglu$ in~\ref{eq:subglue}.

The inductive hypothesis witnesses that
each restriction $h \fcmp \swb_{(\swb \fapp s)}$ is telescopic.
\end{proof}

It follows from this lemma that the function
$\rest \typing (W \extfnc C) \extfnc \submaps$
from \cref{prop:restrglue}.\ref{restrglue:restrext}
restricts to a function
\begin{equation}\label{eq:restlift}
\rest \typing \alg(\ws,a_C) \extfnc \tLab
\end{equation}
that maps an algebra morphism to its family of restrictions
as defined in~\eqref{eq:restr}.
Because of the bijection in \cref{cor:restrglue},
to see this we could have equally well proved
that every coherent family $F$ \famgldt
such that $\rest\fapp (\glu\fapp F) \approx F$
is a telescopic family.
Our choice is justified by the fact that
functions and their restrictions
are simpler objects than families \famgldt and their restrictions.
The fact that every such $F$ is telescopic
is an immediate consequence of the Lemma above
and \cref{prop:restrglue}.\ref{restrglue:rgid}.

Our aim is now to prove that $\glu$ lifts to an inverse of $\rest$.
To do so we need to prove that a telescopic family is coherent.
Then by \cref{prop:restrglue}.\ref{restrglue:rgid} it will be enough to show
that a telescopic family is the restriction of its own glueing.

First we need some technical lemmas about telescopic functions.
We begin by showing that any two functions
which are telescopic over the same $w$ are necessarily equal.
A more general version is in \cref{lem:teluniqtrsp}.

\begin{lmm}\label{lem:teluniqref}
Let $w \typing W$ and $k,k' \typing \isbtree \fapp w \extfnc C $.
Then
\[
\isTel \fapp k \to \isTel \fapp k' \to k \approx k'.
\]
\end{lmm}

\begin{proof}
This is proven by induction on
$T \jueq \wsupd\fapp (w,k)\fapp G\fapp E \typing \isTel \fapp k$
into the type
\[
\prod_{k'} \isTel \fapp k' \to k \approx k'.
\]
Let $G' \coloneqq \slab \fapp w \fapp k' \fapp T'$
be the coherent family \famglsubdt given by the assumption $T'$
that $k'$ is telescopic, as in~\eqref{eq:cslab}.
Since $k$ and $k'$ are the glueing of $G$ and $G'$, respectively,
it is enough to show that
\[
\sglu_w\fapp G \approx \sglu_w\fapp G'.
\]
Using extensionality of $\sglu_w$ from \cref{lem:subglue}.\ref{subglue:ext}
this reduces to show that the two families $G$ and $G'$
\famglsubdt are equal,
namely that for every $s \typing \isbtree \fapp w $
\[
G_s \approx G'_s.
\]
This is the inductive hypothesis~\eqref{eq:telIH}
applied to the right-hand function and the proof from~\eqref{eq:cslabTel}
that it is telescopic.
\end{proof}

\begin{crl}\label{cor:tlabssing}
The setoid $\tLab$ is a subsingleton, that is,
$F \approx F'$ for any two telescopic families $F$ and $F'$.
\end{crl}

\begin{proof}
Straightforward from \cref{def:telstd} and \cref{lem:teluniqref}.
\end{proof}

The next lemma shows that telescopic functions are stable under transport over $W$.

\begin{lmm}\label{lem:teltransp}
Let $\gamma \typing w \approx_W w'$  and $k \typing \isbtree \fapp w \extfnc C$.
Then
\[
\isTel_w \fapp k \to \isTel_{w'} \fapp (k \fcmp \isbtree_{\gamma^{-1}}).
\]
\end{lmm}

\begin{proof}
This is proven by induction on
$T \jueq \wsupd\fapp (w,k)\fapp G\fapp E \typing \isTel \fapp k$
into the type
\[
\prod_{(w' \typing W)} \prod_{(\gamma \typing w \approx w')}%
  \isTel \fapp w' \fapp (k \fcmp \isbtree_{\gamma^{-1}})
\]
and we work towards applying $\wsupd$,
that is conditions~\eqref{tel:lab} and~\eqref{tel:tel}
in \cref{def:tel}.

Since $k$ is obtained applying $\sglu_w$ to $G$ by~\eqref{eq:sglue},
it is
\[
k \fcmp \isbtree_{\gamma^{-1}} \approx \sglu_{w'}\fapp (\csLab_\gamma \fapp G)
\]
by \cref{lem:subglue}.\ref{subglue:coh}.
It remains to provide a branching function to establish condition
\ref{def:tel}.\ref{tel:tel}.
This amounts to show that, for each $s' \typing \isbtree \fapp w'$,
the function
\[
\csLab_\gamma \fapp G_{s'} \jueq G_s \fcmp \isbtree_{(\swbcoh \fapp \gamma^{-1}\fapp s')}
\typing \isbtree\fapp (\swb_{w'}\fapp s') \extfnc C
\]
from~\eqref{eq:cohssmaptrsp} is telescopic, where
$s \coloneqq \isbtree_{\gamma^{-1}} \fapp s' \typing \isbtree\fapp w$.
This is the inductive hypothesis in~\eqref{eq:telIH} applied to
\[
(\swbcoh\fapp \gamma^{-1} \fapp s')^{-1}
\typing \swb_w\fapp s \approx \swb_{w'}\fapp s'
\]
from \cref{rem:cohbr}.
\end{proof}

\begin{lmm}\label{lem:teluniqtrsp}
Let $\gamma \typing w \approx_W w'$, $k \typing \isbtree \fapp w \extfnc C $
and $k' \typing \isbtree \fapp w' \extfnc C$.
Then
\[
\isTel_w \fapp k \to \isTel_{w'} \fapp k' \to k \approx k' \fcmp \isbtree_\gamma.
\]
\end{lmm}

\begin{proof}
Straightforward from \cref{lem:teltransp,lem:teluniqref}.
\end{proof}

Recall that a family \famgldt is telescopic
if each of its components is a telescopic function.
We see in particular that:

\begin{crl}\label{cor:telcoh}
Telescopic families are coherent,
that is, the setoid $\tLab$ is a subsetoid of $\cLab$.
\end{crl}

In the following two results
we make sure that $\glu$ maps telescopic families to algebra morphisms.

\begin{lmm}\label{lem:cmphrisrd}
Let $F$ be a telescopic family.
Then the family $\rest\fapp (\glu\fapp F)$ \famgldt is also telescopic.
\end{lmm}

\begin{proof}
We need to show that, for every $w \typing W$, the function
\[
\rest\fapp (\glu\fapp F)\fapp w \jueq (\glu\fapp F) \fcmp \swb_w
\typing \isbtree\fapp w \extfnc C
\]
is telescopic,
and we do so by proving that it satisfies
conditions~\eqref{tel:lab} and~\eqref{tel:tel}
in \cref{def:tel}.

As family \famglsubdt for $w$ we can take the family $G$
defined by $G_s \coloneqq F_{(\swb_w\fapp s)}$.
Its coherence follows as for \cref{cor:telcoh} from \cref{lem:teluniqtrsp}
and the fact that $F$ is telescopic.
For $s \typing \isbtree\fapp w$ it is
\[\begin{split}
\rest\fapp (\glu\fapp F)\fapp w\fapp s
&\jueq
(\glu\fapp F)\fapp (\swb_w\fapp s)
\\&\jueq
a_C(\wn\fapp (\swb_w\fapp s), F_{(\swb_w\fapp s)})
\\&\jueq
\sglu_w\fapp G
\end{split}\]
by definition of $\glu$ in~\eqref{eq:glue} and
definition of $\sglu$ in~\eqref{eq:subglue}.
Since for every $w$, the function $F_w$ is telescopic,
this is the case in particular for $G_s \jueq F_{\swb_w\fapp s}$
for every $s \typing \isbtree\fapp w$.
\end{proof}

\begin{crl}\label{cor:gluam}
For every telescopic family $F$, it is
\[
\rest\fapp (\glu\fapp F) \approx F
\]
and the function $\glu\fapp F$ is in $\alg(\ws,a_C)$, that is
\[
\glu\fapp F \circ \ws \approx a_C \circ \poly{B}(\glu\fapp F).
\]
\end{crl}

\begin{proof}
The first equality follows from \cref{lem:cmphrisrd}
and \cref{lem:teluniqref}.
So the function $\glu\fapp F$ is an algebra morphism
by \cref{prop:restrglue}.\ref{restrglue:rgid}.
\end{proof}

It follows that the function
$\glu \typing \cLab \extfnc (W \extfnc C)$
defined in~\eqref{eq:glue} restricts to a function
\begin{equation}\label{eq:gluelift}
\glu \typing \tLab \extfnc \alg(\ws,a_C)
\end{equation}
from telescopic families to algebra morphisms.

\begin{thr}\label{thm:commchar}
The functions $\rest$ from~\eqref{eq:restlift} and $\glu$ from~\eqref{eq:gluelift}
\[
\xycenterm[C=5em]{
\alg(\ws,a_C)	\ar@{=>}@<5pt>[r]^-{\rest}	&%
\tLab	\ar@{=>}@<5pt>[l]^-{\glu}
}
\]
establish a bijection between the setoid of algebra morphisms
from $\ws$ to $a_C$ and the setoid of telescopic families \famgldtin{C}.
\end{thr}

\begin{proof}
\Cref{prop:restrglue}.\ref{restrglue:grid} yields that
$\glu\fapp (\rest\fapp h) \approx h$
for every $h \typing \alg(\ws,a_C)$.
The other equality is in \cref{cor:gluam}.
\end{proof}

As an immediate consequence, we see from \cref{cor:tlabssing} that
the setoid of algebra morphisms is a subsingleton.
It is then clear that $\alg(\ws,\ws)$ is a singleton,
and that the only telescopic family in this case
is $\swb$.
It only remains to construct an inhabitant in the general case.

\begin{prp}\label{prop:telglue}
There is a telescopic family \famgldtin{C}, \ie a term
\[
\telglue \typing \tLab.
\]
\end{prp}

\begin{proof}
Given an extensional tree $w \typing W$,
we need to show that  there is a telescopic function
\[
\telglue_w \typing \isbtree \fapp w \extfnc C.
\]
The proof is by $\wty$-elimination on $w \typing \wty[|B|]$ into the type
\[
\wper{B} \fapp w \fapp w	\,\ \longrightarrow%
  \sum_{k \typing \isbtree \fapp w \extfnc C} \isTel \fapp k.
\]
As in the proof of \cref{lem:restristel},
we may assume without loss of generality
that all the trees in the proof are extensional.

We have a tree of the form $w \jueq \wsup \fapp a \fapp f$
and the inductive hypothesis consists,
for every $s \typing \isbtree\fapp w$, of a telescopic function
\[
G_s \typing \isbtree\fapp (\swb_w\fapp s) \extfnc C.
\]
This is the same as a telescopic family $G$ \famglsubdt for $w$,
see~\cref{rem:telalt}.

It follows from \cref{lem:teluniqtrsp} that $G$ is coherent,
so we may define
\[
\telglue _w \coloneqq \sglu_w\fapp G \typing \isbtree \fapp w \extfnc C
\]
as the glueing of a telescopic family $G$ \famglsubdt for $w$.
It is a telescopic function by definition.
\end{proof}

\subsection{Initiality}
\label{ssec:init}

We have finally reached our main result.

\begin{thr}\label{thm:init}
For every setoid family $B$ over a setoid $A$,
the associated polynomial endofunctor $\poly{B}$ has an initial algebra.

It follows that the category of setoids \catstdt
has initial algebras for polynomial endofunctors.
\end{thr}

\begin{proof}
Given a setoid family $B$ over $A$,
the initial algebra for the polynomial functor $\poly{B}$
is the algebra $\ws_B \typing \poly{B}W \extfnc W$
of extensional trees constructed in \cref{ssec:exttree}.
It follows by \cref{thm:commchar} that,
for every algebra $a_C$,
the only morphism of algebras from $\ws$ to $a_C$
is the function $\glu (\telglue) \typing W \extfnc C$,
where $\glu \typing \tLab \extfnc \alg(\ws_B,a_C)$
is the glueing function from \cref{thm:commchar}
and $\telglue$ is the only telescopic family \famgldtin{C},
that we constructed in \cref{prop:telglue}
and proved unique in \cref{cor:tlabssing}.

It is straightforward to verify that
every polynomial functor $P$ on \catstdt has initial algebra
$\ws_B \fcmp i_W \typing PW \extfnc W$,
where $i$ is a natural isomorphism from
$P$ to $\poly{B}$ as in \cref{def:poly}.
\end{proof}

Gambino and Hyland have proved that,
in a locally cartesian closed category,
dependent W-types can be constructed from non-dependent ones,
thus providing justification for an analogous result in
extensional type theory claimed in~\cite{PeterssonSynek1989}.
As a consequence of \cref{thm:init},
we expect that setoids have dependent W-types too,
\ie that the category of \catstdt has
initial algebras for \emph{dependent} polynomial endofunctors.

\section{Trees on discrete setoids}
\label{sec:disctree}

In this section we work with the additional assumption that
our theory has identity types~\cite[Chapter~20.4]{NPS1990},
see also~\cite[Chapter~1.12]{HoTTbook}.
We unfold the definition of the setoid of extensional trees
in the case of a \emph{discrete} setoid family,
\ie a setoid family coming from
an extensional function between discrete setoids.
We then compare it to the discrete setoid
on the W-type of the underlying type family.
The main result is \cref{thm:freestd}.

In comparing the two, we will also make use of function extensionality.
This is because the functions of immediate subtrees of two matching trees
are only required to be pointwise equal,
as in the right-hand side of~\eqref{condW2}.
Therefore, rather than to compare the two directly,
we prefer to identify another (total) equivalence relation on $\wty[|B|]$,
which is logically equivalent to the identity type on $\wty[|B|]$
in the presence of function extensionality,
and to compare it to $\wper{B}$ in the general case.

Recall that $x =_X x'$ denotes
the identity type of two terms $x,x' \typing X$.
For a type family $Y \typing X \to \ttunivm$,
we denote as $Y_{\xi} \typing Y\fapp x \to Y\fapp x'$ the transport function term
given by identity elimination on $\xi \typing x =_X x'$,
and we write $y =_{\xi} y'$ for the identity type
$Y_{\xi}\fapp y =_{(Y\fapp x')} y'$.

\subsection{Discrete setoid families}

We begin by identifying those setoid families that
correspond, under \cref{thm:stdfam}, to functions between discrete setoids.

\begin{dfn}\label{def:freestdfam}
A setoid family $B \typing \sfam X$ over the discrete setoid on $X$
is \emph{discrete over $X$} if, for every $x \typing X$ and $b,b' \typing B\fapp x$,
\[
b \approx_{B\fapp x} b'
\ \longleftrightarrow \
\sum_{l \typing x =_X x} |B|_l\fapp b =_{B\fapp x} b'
\]
where $|B|_{\xi} \typing B\fapp x \to B\fapp x$ is the transport function
of the underlying type family $|B|$ of $B$
given by elimination of the identity type.
\end{dfn}

The equivalence in \cref{thm:stdfam} restricts to an equivalence
between the full subcategory of $\catstdm/X$ on the functions into $X$
from a discrete setoid,
and the full subcategory of $\sfam X$ on the discrete setoid families over $X$.

It follows in particular that every type family
$Y \typing X \to \ttunivm$ gives rise to a discrete setoid family over $X$,
by equipping each fibre $Y\fapp x$ with the equality
\[
y \approx y' \coloneqq \sum_{\xi \typing x =_X x} Y_{\xi}\fapp y = y'
\]
where $Y_{\xi} \typing Y\fapp x \to Y\fapp x$ is
the transport function given by elimination of the identity type.

\subsection{Pointwise equality of trees}
Next, we look for an equivalent characterisation of
the discrete setoid on a W-type in the presence of
function extensionality.

\begin{dfn}\label{def:eqw}
Let $Y \typing X \to \ttunivm$ be a family of small types.
Two trees $w \jueq \wsup\fapp x\fapp f$ and $w' \jueq \wsup\fapp x'\fapp f'$ 
in $\wty[Y]$ are \emph{pointwise equal} if
$\xi \typing x =_X x'$ and, for every $y \in Y\fapp x$ and $y' \typing Y\fapp x'$,
if $y =_{\xi} y'$ then the immediate subtrees
$f\fapp y$ and $f'\fapp y'$ are pointwise equal.

The type family $\eqWty_Y \typing \wty[Y] \to \wty[Y] \to \ttunivm$
of proofs that $w$ and $w'$ are pointwise is defined
as the (curried version of the) dependent W-type on
$I \coloneqq \wty[Y] \times \wty[Y]$
with family of nodes
\[
N \coloneqq \lambda (w,w'). \wn\fapp w =_X \wn\fapp w'
\typing \wty[Y] \times \wty[Y] \to \ttunivm
\]
branching family
\[
Br \coloneqq \lambda (w,w'),p. \sum_{y,y'} y =_{\xi} y'
\typing \prod_I N\fapp (w,w') \to \ttunivm
\]
and arity function
\[
ar := \lambda (w,w'),p,(y,y',q). (\wb_w\fapp y, \wb_{w'}\fapp y')
\typing \prod_I \prod_N Br((w,w'),p)
\to \ttunivm.
\]
\end{dfn}

It is possible to show that $\eqWty_Y$ is a symmetric and transitive relation on $\wty[Y]$
as for $\wper{B}$ in \cref{prop:wpersymtra}.
Reflexivity follows from the following lemma.

\begin{lmm}\label{lem:eqwenc}
For every $w,w' \typing \wty[Y]$
\[
w =_{\wty[Y]} w' \longrightarrow \eqWty_Y\fapp w\fapp w'.
\]
\end{lmm}

Thus we have a setoid $\widetilde{\wty}_Y$ of pointwise-equal trees,
that is the type $\wty[Y]$ with $\eqWty$ as equality,
together with an extensional function
\begin{equation}\label{eq:eqw}
q_Y \typing \wty[Y] \extfnc \widetilde{\wty}_Y
\end{equation}
from the discrete setoid on $\wty[Y]$
which is in fact a (canonical) quotient map in \catstdt, cf.\ \cref{rem:imgfact}.

\begin{proof}[Proof of \cref{lem:eqwenc}]
The proof is by induction on $w \jueq \wsup\fapp x\fapp f$ and
the inductive hypothesis tells us that if a tree is equal
to an immediate subtree of $w$, then it is pointwise equal to it.

By induction on $r \typing w = w'$, it is enough to show that $\eqWty_Y\fapp w\fapp w$.
A canonical term is given taking the canonical proof of $x =_X x$ and
proving that, for every $y,y' \typing Y\fapp x$, if $y =_{(Y\fapp x)} y'$
then the trees $f\fapp y$ and $f\fapp y'$ are pointwise equal.
This follows from the inductive hypothesis using the fact that every function term---%
$f$ in particular, is extensional \wrt the identity type.
\end{proof}

Unsurprisingly,
it is for the converse that we need function extensionality.

\begin{lmm}\label{lem:eqwdec}
Assume that the theory has function extensionality,
\ie that there is a constant $\mathsf{funext}$ as given in the rule below.
\begin{equation}\label{eq:funext}
\begin{prooftree}
X,Y \typing \ttunivm
\hspace{2em}
f,g \typing X \to Y
\hspace{2em}
H \typing \prod_{x \typing X} f\fapp x =_Y g\fapp x
\justifies \rule{0pt}{10pt}
\mathsf{funext}\fapp H \typing 
f =_{(X \to Y)} g
\hspace{2em}
\end{prooftree}
\end{equation}
Then for every $w,w' \typing \wty[Y]$
\[
\eqWty_Y\fapp w\fapp w' \longrightarrow w =_{\wty[Y]} w'.
\]
\end{lmm}

\begin{proof}
Let $w \jueq \wsup\fapp x\fapp f$ and $w' \jueq \wsup\fapp x'\fapp f$'.
The proof is by induction on the term $E \typing \eqWty\fapp w\fapp w'$.
To conclude $w = w'$ it is enough to show that there is $\xi \typing x =_X x'$
such that $f = f' \fcmp Y_{\xi}$.
The existence of $p$ follows by definition of $\eqWty\fapp w\fapp w'$.
By function extensionality and the inductive hypothesis it is then enough to show that,
for every $y \typing Y\fapp x$,
the immediate subtrees
$f\fapp y$ and $f' \fapp (Y_{\xi}\fapp y)$ are pointwise equal.
This also follows by definition of $\eqWty\fapp w\fapp w'$.
\end{proof}

It follows that, in the presence of function extensionality,
the discrete setoid on $\wty[Y]$ is isomorphic to the setoid
$\widetilde{\wty}_Y$ of pointwise-equal trees
via the function $q_Y$ in~\eqref{eq:eqw}.

\subsection{Equivariant trees}

Now we look for a characterisation of the setoid $W_B$ of extensional trees
of a discrete setoid family $B$ as subsetoid of $\widetilde{W}_{|B|}$,
\ie as a setoid of pointwise equal trees.

\begin{dfn}\label{def:eqv}
Let $Y \typing X \to \ttunivm$ be a type family.
A tree $w \jueq \wsup\fapp x\fapp f \typing \wty[Y]$ is \emph{equivariant}
if, for every $l \typing x =_X x$ and every $y \typing Y\fapp x$,
the immediate subtrees $f\fapp (Y_l\fapp y)$ and $f\fapp y$ are pointwise equal,
and every immediate subtree of $w$ is equivariant.

The type family $\isEqv \typing \wty[Y] \to \ttunivm$ of proofs
that a tree is equivariant is defined as the
dependent W-type in $I \coloneqq \wty[Y]$ with families of nodes and branches
\[
N \coloneqq \lambda w. \prod_{(l \typing \wn\fapp w =_X \wn\fapp w)}
\prod_{(y \typing Y\fapp (\wn\fapp w))} 
\eqWty_Y (f\fapp (Y_l\fapp y), f\fapp y),
\hspace{3em}
Br \coloneqq \lambda w,\_. Y\fapp (\wn\fapp w)
\]
and arity function
\[
ar \coloneqq \lambda w,\_,y. \wb_w\fapp y.
\]

Let $\Eqv_Y$ be the subsetoid of $\widetilde{\wty}_Y$ on the equivariant trees,
that is to say, two equivariant trees are equal in $\Eqv_Y$
if they are pointwise equal.

Let also $\Eqv^=_Y$ be the subsetoid of $\wty[Y]$ on the equivariant trees,
that is,
two equivariant trees $w$ and $w'$ ar equal in $\Eqv^=_Y$
if $w =_{\wty[Y]} w'$.
\end{dfn}

Let $B \typing \sfam X$ be a discrete setoid family
over the discrete setoid on $X \typing \ttunivm$.
Note first that, by its very definition,
for every $\xi \typing x =_X x'$ and every $b \typing B\fapp x$,
there is $l' \typing x' =_X x'$ such that
\begin{equation}\label{eq:freetransp}
|B|_{l'}\fapp (B_{\xi}\fapp b) =_{B\fapp x'} |B|_{\xi}\fapp b.
\end{equation}
Recall that we write $B_{\xi}$ for the transport of the setoid family $B$,
and $|B|_{\xi}$ for the transport of the underlying type family $|B|$---%
the former being given by definition and the latter
by elimination of the identity type.
The matching relation $\wper{B}$ was defined in \ref{def:Wper}.

\begin{lmm}\label{lem:wper2eqv}
Let $B \typing \sfam X$ be a discrete setoid family
over the discrete setoid on $X \typing \ttunivm$.
\begin{enumerate}
\item
For every $w,w' \typing \wty[|B|]$,
\[
\wper{B}\fapp w\fapp w' \longrightarrow \eqWty_{|B|}\fapp w\fapp w'.
\]
\item
Every extensional tree is equivariant, \ie for every $w \typing \wty[|B|]$,
\[
\wper{B}\fapp w\fapp w \longrightarrow \isEqv\fapp w.
\]
\end{enumerate}
\end{lmm}

\begin{proof}
(1)
The proof is by induction on
$w \jueq \wsup\fapp x\fapp f \typing \wty[|B|]$ and the
inductive hypothesis tells us that if a tree is matching with an immediate
subtree $f\fapp b$ of $w$, then it is pointwise equal to $f\fapp b$.

Let then $w' \jueq \wsup\fapp x'\fapp f'$ be a tree an suppose that
$w$ and $w'$ are matching.
To prove that $w$ and $w'$ are pointwise equal we need first to provide a proof
$\xi \typing x =_X x'$, which is given by $\alpha$ in \cref{rem:charWper}.
We then need to show that for every $b \typing B\fapp x$ and $b' \typing B\fapp x'$,
if $B_{\xi}\fapp b \approx_{B\fapp x'} b'$, then the immediate subtrees $f \fapp b$
and $f'\fapp b'$ are pointwise equal.
By the inductive hypothesis, it is enough to show that the trees $f \fapp b$
and $f'\fapp b'$ are matching, which holds by definition.

(2)
The proof is again by induction on $w \typing \wty[|B|]$, into the type
\[
\wper{B}\fapp w\fapp w \longrightarrow \isEqv\fapp w.
\]
By the inductive hypothesis and the fact that every subtree of an extensional subtree
is extensional~\ref{lem:exttree},
it follows that every immediate subtree of $w$ is equivariant.
To prove that $w$ is equivariant,
it only remains to show that for every $l \typing x =_X x$
and every $b\typing B\fapp x$,
the subtrees $f\fapp (|B|_l\fapp b)$ and $f\fapp b$ are pointwise equal.
By (1) just proved and the fact that
$f$ is extensional~\ref{lem:exttree},
it is enough to prove that $|B|_l\fapp b \approx_{B\fapp x} b$.
This is immediate by \cref{def:freestdfam}
and standard properties of transport along identity proofs.
\end{proof}

It follows from \cref{lem:wper2eqv} that there is an extensional function
\begin{equation}\label{eq:wper2eqv}
j_B \typing W_B \extfnc \Eqv_{|B|}
\end{equation}
which is the identity on the underlying trees.
The next lemma proves that $j_B$ is injective.

\begin{lmm}\label{lem:eqv2wper}
Let $B \typing \sfam X$ be a discrete setoid family
over the discrete setoid on $X \typing \ttunivm$,
and let $w,w' \typing \wty[|B|]$ be equivariant trees.
Then
\[
\eqWty_{|B|}\fapp w\fapp w' \longrightarrow \wper{B}\fapp w\fapp w'.
\]
\end{lmm}

\begin{proof}
This is proven by induction on the proof $E$ that
$w \jueq \wsup\fapp x\fapp f$ is equivariant,
into the type
\[
\prod_{w' \typing \wty[Y]} \isEqv\fapp w'
	\to \eqWty_{|B|}\fapp w\fapp w' \to \wper{B}\fapp w\fapp w'.
\]
Let then $w' \jueq \wsup\fapp x'\fapp f'$ be equivariant and suppose that
$w$ and $w'$ are pointwise equal.
In particular, there is $\xi \typing x =_X x'$.
So to prove that $w$ and $w'$ are matching it is enough to show that
for every $b \typing B\fapp x$ and $b' \typing B\fapp x'$,
if $B_{\xi}\fapp b \approx_{B\fapp x'} b'$
then the immediate subtrees $f\fapp b$ and $f'\fapp b'$ are matching.
This follows by the inductive hypothesis once we know that
$f'\fapp b'$ is equivariant, and that $f\fapp b$ and $f'\fapp b'$
are pointwise equal.

The former condition follows by \cref{def:eqv}.
For the latter we have
\[\begin{split}
f\fapp b &\approx_{\eqWty} f' (|B|_{\xi}\fapp b)
\\
&\approx_{\eqWty}
f'\fapp (B_{\xi}\fapp b)
\\
&\approx_{\eqWty}
f'\fapp b'
\end{split}\]
using, in order,
the fact that subtrees of $w$ and $w'$ are pointwise equal by \cref{def:eqw};
condition \eqref{eq:freetransp}, \cref{lem:eqwenc}
and equivariance of $w'$;
and \cref{def:freestdfam},
\cref{lem:eqwenc} and equivariance of $w'$.
\end{proof}

\begin{thr}\label{thm:freestd}
Let $B \typing \sfam X$ be a discrete setoid family
over the discrete setoid on $X$.
\begin{enumerate}[beginpenalty=99,midpenalty=99]
\item
The function $j_B$ in \eqref{eq:wper2eqv} is a bijection between
the setoid of extensional trees $W_B$
and the setoid of equivariant trees $\Eqv_{|B|}$
on the underlying type family $|B|$.
\item
Assuming function extensionality~\eqref{eq:funext},
the setoid of extensional trees $W_B$ on $B$ is in bijection
with the subsetoid  $\Eqv^=_Y$ of the
(discrete setoid on the) W-type $\wty[|B|]$.
\end{enumerate}
\end{thr}

\begin{proof}
(1)
This is immediate from \cref{lem:wper2eqv,lem:eqv2wper}.

(2)
By \cref{lem:eqwenc,lem:eqwdec}, the discrete setoid of trees $\wty[|B|]$
is isomorphic to the setoid $\widetilde{\wty}_{|B|}$ of pointwise-equal trees.
Now the claim follows by (1)
since the setoid of equivariant trees $\Eqv_{|B|}$
is a subsetoid of $\widetilde{\wty}_{|B|}$ by definition.
\end{proof}

Even if assuming function extensionality does not seem to ensure, in general,
that the setoid $W_B$ of extensional trees
\wrt a discrete setoid family $B$ is discrete,
this happens whenever the base type is a 0-type,
\ie a type with decidable equality.

\begin{crl}\label{cor:hset}
Let $B \typing \sfam X$ be a discrete setoid family
over the discrete setoid on $X$,
and suppose that the type $X$ is a 0-type.
\begin{enumerate}
\item
Every tree in $\wty[|B|]$ is extensional \wrt $B$.
In particular, the setoid $W_B$ of extensional trees is in bijection
with the setoid $\widetilde{\wty}_{|B|}$ of pointwise-equal trees.
\item
Assuming function extensionality~\eqref{eq:funext},
the setoid $W_B$ of extensional trees is in bijection with the
discrete setoid on the W-type $\wty[|B|]$.
\end{enumerate}
\end{crl}

\begin{proof}
We only need to prove that every tree is extensional,
as all the other claims follow from \cref{thm:freestd}.

A tree is extensional \iff it is equivariant by \cref{lem:wper2eqv,lem:eqv2wper},
and it is straightforward to prove that every tree in $\wty[|B|]$ is
equivariant, by induction on the tree and using the assumption that $X$ is a 0-type.
\end{proof}

\section*{Acknowledgements}

The work described in this paper would have not been possible without the support of
my supervisor Erik Palmgren and, in particular,
his Coq library on setoids and setoid families.
The main result of this paper was presented at the
Workshop on Types, Homotopy Type Theory, and Verification,
held at the Hausdorff Research Institute for Mathematics in Bonn in June 2018,
and a first version was completed
while I was hosted at the same Institute in July 2018.
I thank the organisers of the workshop for giving me the opportunity to speak,
the participants for valuable feedback
and the Institute for excellent working conditions.
Financial support from the Royal Swedish Academy of Sciences and
the K\&A Wallenberg Foundation is also acknowledged.
I am grateful to Peter Dybjer for bibliographic advice
and to the anonymous referees for extremely useful comments.
Prooftrees were typeset using Paul Taylor's macros package.

\bibliographystyle{alpha}
\bibliography{w-types-std-references}

\end{document}

%% file: w-types-std-macros.tex
\theoremstyle{plain}
\newtheorem{thr}[thm]{Theorem}
\newtheorem{prp}[thm]{Proposition}
\newtheorem{lmm}[thm]{Lemma}
\newtheorem{crl}[thm]{Corollary}

\theoremstyle{definition}
\newtheorem{dfn}[thm]{Definition}
\newtheorem{rmk}[thm]{Remark}
\newtheorem{xmps}[thm]{Examples}

\crefname{thr}{Theorem}{Theorems}
\crefname{prp}{Proposition}{Propositions}
\crefname{lmm}{Lemma}{Lemmas}
\crefname{crl}{Corollary}{Corollaries}
\crefname{dfn}{Definition}{Definitions}
\crefname{rmk}{Remark}{Remarks}
\crefname{xmps}{Examples}{Examples}

\newcommand{\ie}{\textit{i.e.}\xspace}
\newcommand{\eg}{\textit{e.g.}\xspace}
\renewcommand{\iff}{if and only if\xspace}
\newcommand{\wrt}{with respect to\xspace}
\newcommand{\mltt}{Martin-L\"of type theory\xspace}
\newcommand{\epsi}{\varepsilon}

\newcommand{\abscatm}[1]{\mathbb{#1}}		
\newcommand{\catam}{\abscatm{A}}		\newcommand{\catat}{$\catam$\xspace}
\newcommand{\catbm}{\abscatm{B}}		\newcommand{\catbt}{$\catbm$\xspace}
\newcommand{\catcm}{\abscatm{C}}		\newcommand{\catct}{$\catcm$\xspace}
\newcommand{\conccatm}[1]{\mathsf{#1}}		
		
\newcommand{\catstdm}{\conccatm{Std}}		\newcommand{\catstdt}{$\catstdm$\xspace}
\newcommand{\xycenterm}[2][=2em]{\vcenter{\hbox{\xymatrix@#1{#2}}}}

\newcommand{\typing}{:}
\newcommand{\ttunivm}{\mathsf{U}}	\newcommand{\ttunivt}{$\ttunivm$\xspace}
\newcommand{\proj}{\textup{pr}}
\newcommand{\ttemptym}{\mathsf{0}}	\newcommand{\ttemptyt}{$\ttemptym$\xspace}
\newcommand{\ttunitm}{\mathsf{1}}	\newcommand{\ttunitt}{$\ttunitm$\xspace}
\newcommand{\wty}[1][]{\mathsf{W}_{#1}}
\newcommand{\dwty}[1][]{\mathsf{DW}_{#1}}
\newcommand{\wsup}{\mathsf{sup}}
\newcommand{\wsupd}{\mathsf{dsup}}
\newcommand{\wn}{\mathsf{node}}
\newcommand{\wb}{\mathsf{ist}}
\newcommand{\ws}{\mathsf{s}}
\newcommand{\wus}{\mathsf{us}}
\newcommand{\dwn}{\mathsf{dnode}}
\newcommand{\dwb}{\mathsf{dist}}
\newcommand{\wbcoh}{\mathsf{\wb\text{-}coh}}
\newcommand{\swbcoh}{\mathsf{\swb\text{-}coh}}
\newcommand{\swb}{\mathsf{m}}
\newcommand{\istq}{\mathsf{e}}
\newcommand{\extfnc}{\Rightarrow}
\newcommand{\extpf}{\mathop{\Rsh}}
\newcommand{\wrec}[1][]{\mathsf{rec}^{\mathsf{W}}_{#1}}
\newcommand{\dwrec}[1][]{\mathsf{rec}^{\mathsf{DW}}_{#1}}
\newcommand{\jueq}{\equiv}

\newcommand{\fapp}{\,}
\newcommand{\fcmp}{\circ}

\newcommand{\sfam}{\mathsf{Fam}\,}

\newcommand{\isCoh}[1]{\mathsf{isCoh}_{#1}}

\newcommand{\eqWty}{\mathsf{ptEq}}
\newcommand{\isEqv}{\mathsf{isEquivariant}}
\newcommand{\Eqv}{\mathsf{EqvTrees}}
\newcommand{\isbtree}{\mathsf{ISTrees}}

\newcommand{\csLab}{\mathsf{CohGlueSubD}}
\newcommand{\submaps}{\mathsf{GlueD}}
\newcommand{\cLab}{\mathsf{CohGlueD}}
\newcommand{\sglucoh}{\mathsf{subglue\text{-}coh}}
\newcommand{\isTel}{\mathsf{isTelescopic}}
\newcommand{\telglue}{\mathsf{telfam}}

\newcommand{\slab}{\mathsf{csfam}}
\newcommand{\telsglu}{\mathsf{isglue}}
\newcommand{\telslab}{\mathsf{telsfam}}
\newcommand{\alg}{\mathsf{Alg}}
\newcommand{\tLab}{\mathsf{TelGlueD}}
\newcommand{\rest}{\mathsf{restr}}
\newcommand{\glu}{\mathsf{glue}}
\newcommand{\sglu}{\mathsf{subglue}}
\DeclareMathOperator{\Obj}{Obj}
\DeclareMathOperator{\Hom}{Hom}
\newcommand{\efun}[1]{\mathsf{Fun}(#1)}
\newcommand{\poly}[1]{\mathsf{P}}
\newcommand{\polysub}[1]{\mathsf{P}_{\!\!#1}}
\newcommand{\wper}[1]{\mathsf{W}^{\mathsf{per}}_{#1}}
\newcommand{\famgldt}{of glueing data\xspace}
\newcommand{\famgldtin}[1]{of glueing data in $#1$\xspace}
\newcommand{\famglsubdt}{of glueing subdata\xspace}

\newcommand{\famtel}{telescopic family\xspace}
\newcommand{\famstel}{telescopic families\xspace}